\RequirePackage{fix-cm}

\documentclass[smallextended]{svjour3}       % onecolumn (second format)
\smartqed  % flush right qed marks, e.g. at end of proof
\usepackage{graphicx}
\usepackage{amsmath,amssymb,amsfonts}%
\usepackage{bm}
\usepackage{hyperref}
\usepackage{mathabx}
\usepackage{subcaption}
\usepackage{xcolor}

\renewcommand{\bar}[1]{\widebar{#1}} 
\renewcommand{\hat}[1]{\widehat{#1}} 
\renewcommand{\Tilde}[1]{\widetilde{#1}} 
\newcommand{\LRp}[1]{\left( #1 \right)} % parens
\newcommand{\LRs}[1]{\left[ #1 \right]} % square brackets
\newcommand{\LRv}[1]{\lvert #1 \rvert}
\newcommand{\rmd}{{\rm{d}}}
\newcommand{\rmt}{{\rm{t}}}
\newcommand{\p}{\partial}
\newcommand{\f}{\frac}
\newcommand{\csp}{,\qquad}
\begin{document}

\title{Model Order Reduction Techniques for the Stochastic Finite Volume Method
% \thanks{Grants or other notes
% about the article that should go on the front page should be
% placed here. General acknowledgments should be placed at the end of the article.}
}

\author{Ray Qu \and
        Jesse Chan \and
        Svetlana Tokareva%etc.
}

% \authorrunning{Short form of author list} % if too long for running head

\institute{R. Qu \at
              Department of Computational Applied Mathematics and Operations
Research, Rice University, 6100 Main St MS-134, Houston, TX 77005, USA \\
              \email{Ray.Qu@rice.edu}           
           \and
           J. Chan \at
              Oden Institute for Computational Engineering and Sciences, The University of Texas at Austin, 201 E 24th St POB 4.102, Austin, TX 78712, USA \\
            \email{jesse.chan@oden.utexas.edu}
            \and
            S. Tokareva \at
            Theoretical Division, Los Alamos National Laboratory, MS B284, Los
Alamos, NM 87545, USA \\
\email{tokareva@lanl.gov}
}

\date{Received: date / Accepted: date}

\maketitle

\begin{abstract}
The stochastic finite volume method (SFV method) is a high-order accurate method for uncertainty quantification (UQ) in hyperbolic conservation laws. However, the computational cost of SFV method increases for high-dimensional stochastic parameter spaces due to the curse of dimensionality. To address this challenge, we incorporate interpolation-based reduced order model (ROM) techniques that reduce the cost of computing stochastic integrals in the SFV method. Further efficiency gains are achieved through hyper-reduction with a QR factorization-based discrete empirical interpolation method (Q-DEIM). Numerical experiments suggest that this approach can lower both computational cost and memory requirements for high-dimensional stochastic parameter spaces.

\end{abstract}

\section{Introduction}

Hyperbolic systems of partial differential equations are ubiquitous in physics and engineering, where mathematical models are represented by a set of transport dominated (hyperbolic) balance laws. Many efficient numerical methods have been developed to approximate the solutions of systems of conservation laws \cite{LeVeque,Godlewski-Raviart}, e.g. finite volume schemes \cite{morton2007finite} or discontinuous Galerkin methods \cite{krivodonova2004shock,Cockburn2001,HesthavenWarburton2007}. The classical assumption in designing efficient numerical methods for hyperbolic systems is that the initial data, boundary conditions, and coefficients of the model are known exactly, i.e., they are deterministic. However, in many practical applications it is not always possible to obtain exact data due to, for example, measurement, prediction, or modeling errors. 

It is well known that solutions to hyperbolic PDEs may develop discontinuities in finite time even from smooth initial data \cite{Toro,Dafermos2016}. In the context of stochastic hyperbolic PDEs (SPDEs), this leads to the propagation of discontinuities in both physical and stochastic dimensions, posing significant challenges on the numerical approximation of SPDEs.

Numerical techniques have been developed to quantify uncertainty by computing the mean flow and its statistical moments.  They are based on approaches such as Monte Carlo (MC) \cite{MishraSchwab2010}, multilevel Monte Carlo (MLMC) \cite{MishraSchwabSukys2011}, generalized polynomial chaos (gPC) \cite{Poette,Lin-Su-Karniadakis-1}, the probabilistic collocation method (PCM) \cite{Lin-Su-Karniadakis-2}, Godunov schemes \cite{Troyen-2}, multiresolution methods \cite{Kolb2023} and stochastic Galerkin (sG) projections \cite{Troyen-1,Gottlieb-Xiu}.

Uncertainty quantification methods can be roughly classified into intrusive and non-intrusive schemes. \textit{Non-intrusive} UQ methods allow reuse of an existing deterministic code as a black box without any modifications \cite{MishraSchwab2010,MishraSchwabSukys2011,Lin-Su-Karniadakis-2,REAGANA2003,Berveiller2006,HosderWalters}. An example of such an approach is the Monte Carlo method \cite{MishraSchwab2010,MishraSchwabSukys2011}: given a large number of realizations of the random parameters, we generate the corresponding outputs from our deterministic code and then process these data to obtain the statistical mean and variances. The possibility to reuse existing deterministic codes is a clear advantage of non-intrusive approaches. However, methods such as Monte Carlo tend to require a prohibitively large number of evaluations of the deterministic solutions. In contrast, \textit{intrusive} UQ methods require modifications to the algorithms and their implementations in the deterministic computational scheme \cite{Poette,Lin-Su-Karniadakis-1,Troyen-2,Kolb2023,Troyen-1,Gottlieb-Xiu,LeMaitre2010}. For example, a well-known and widely used sG method for hyperbolic systems transforms the original PDEs for the primary variables into a set of PDEs that are defined with respect to the polynomial expansion coefficients \cite{Troyen-1}. This reduces the simulation time but might pose mathematical challenges, such as the loss of hyperbolicity of the transformed PDEs.

Recently, we have developed a semi-intrusive Stochastic Finite Volume (SFV) method \cite{Tokareva14} to quantify uncertainties that arise due to random model parameters in the underlying hyperbolic PDE system, including initial conditions, uncertain constants, and complex time-dependent distributions on boundary conditions.  The SFV method requires some modifications of the deterministic code that is used for solving a standard initial boundary value problem (IBVP) for hyperbolic conservation laws, which only involve additional integration of the numerical fluxes over the cells in the stochastic space and can therefore be considered mild.  This approach preserves the hyperbolicity of the model and, at the same time, is more computationally efficient than the Monte Carlo method.

The main bottleneck in applying the SFV method to high-dimensional SPDEs is \emph{curse of dimensionality:} an exponential growth of algorithmic complexity and memory requirements with an increasing number of stochastic dimensions. Several promising approaches have been proposed to improve the scalability of SFV: an adaptive SFV \cite{Harmon2025} and low-rank tensor-train SFV \cite{WALTON2025}. 
    
In this paper, we explore the interpolation-based reduced order model (ROM) techniques that reduce the cost of computing stochastic integrals in the SFV method. When constructing the reduced basis, we highlight a non-intrusive approach that generates stochastic snapshots and applies weighted essentially non-oscillatory (WENO) reconstruction to these snapshots. Further efficiency gains from hyper-reduction are achieved through the discrete empirical interpolation method (DEIM) \cite{Chaturantabut10} with node selection via the QR factorization with column pivoting, known as Q-DEIM \cite{Drmac16}.

This paper is organized as follows: Section \ref{sec:SFV method} reviews the SFV method for conservation laws in a one-dimensional physical space. Section \ref{sec:WENO} presents and compares two reconstruction strategies within the framework of the SFV method. Section \ref{sec:reduced} applies an interpolation-based reduced-order modeling technique. Section \ref{sec:num_exp} provides numerical experiments for Burgers' and compressible Euler equations that assess the performance of the proposed approach.

\section{Stochastic finite volume method}
\label{sec:SFV method}
Our work considers the stochastic hyperbolic systems of conservation laws under a one-dimensional spatial domain and $q$-dimensional parametrized probability space, which is given by the parametric form:

\begin{equation}
\begin{array}{c}
\partial_t\bm{u} + \p_x \bm{F}\LRp{\bm{u},\bm{y}}  = \bm{0}\csp x\in D_x \subset \mathbb{R} \csp \bm{y}\in D_y\subset \mathbb{R}^q\csp t>0; \\
\bm{u}\LRp{x,0,\bm{y}} = \bm{u}_0\LRp{x,\bm{y}}\csp x\in D_x\subset \mathbb{R} \csp \bm{y}\in D_y\subset \mathbb{R}^q,
\end{array}
 \label{eq:stochastic_conservation_law}
\end{equation}
where $\bm{u}=\LRs{u_1,...,u_n}^T,$ $\bm{F}  = \LRs{f_1,...,f_n}^T.$ 

Let $\cup_{i=1}^{N_x}K_x^i$ be the partition of the computational domain $D_x$, each with length $\LRv{K_x^i}$, and $\cup_{j=1}^{N_y}K_y^j$ be the Cartesian grid in the domain $D_y$. We further assume the existence of the probability density function $\mu(\bm{y})$ and compute the expectation of the $p$-th solution component of \eqref{eq:stochastic_conservation_law}
as 
\[
\mathbb{E}\LRs{u_p} = \int_{D_y}u_p\mu(\bm{y})\: d\bm{y}\csp p = 1,...,n.
\]

The SFV method can be obtained by integrating \eqref{eq:stochastic_conservation_law} over each physical-stochastic control volume $K_x^i \times K_y^j$:
\[
\int_{K_y^j}\int_{K_x^i}\partial_t \bm{u}\mu(\bm{y})\:dxd\bm{y} + \int_{K_y^j}\int_{K_x^i} \p_x\bm{F}\LRp{\bm{u},\bm{y}}\mu(\bm{y})\:dxd\bm{y} = \bm{0}. 
\]
Introducing the cell average and stochastic cell measure
\[
\bm{U}_{i,j}(t) = \f{1}{\LRv{K_x^i}\LRv{K_y^j}}\int_{K_y^j}\int_{K_x^i} \bm{u}\LRp{x,t,\bm{y}}\mu(\bm{y})\:dxd\bm{y}\csp \LRv{K_y^j} = \int_{K_y^j}\mu(\bm{y})\:d\bm{y}
\]
and performing the partial integration over $K_x^i$, we get
\[
    \f{d \bm{U}_{i,j}}{dt} + \f{1}{\LRv{K_x^i}\LRv{K_y^j}}\int_{K_y^j}\LRs{\int_{\partial K_x^i} \bm{F}\LRp{\bm{u}}\cdot \bm{n} \: dS}\mu(\bm{y}) \: d\bm{y} = \bm{0}.
\]
Next, we approximate the flux over cell boundaries $\bm{F}\LRp{\bm{u}}\cdot \bm{n}$ with a numerical flux
$$ \hat{\bm{F}}\LRp{\Tilde{\bm{u}}_L\LRp{x,t,\bm{y}}, \Tilde{\bm{u}}_R\LRp{x,t,\bm{y}}},$$ where $\Tilde{\bm{u}}_{L,R}$ denote the reconstructed solution values at the edges of $K_x^i$. For simplicity of notation, we omit the time $t$ for the remainder of this work. We further denote $\hat{\bm{F}}_{i\pm\f{1}{2}}$ as the values of $\hat{\bm{F}}$ at the left $\LRp{i-\f{1}{2}}$ and right $\LRp{i+\f{1}{2}}$ edges of each spatial cell $K_x^i$. The numerical flux integral can then be simplified into the following integral over stochastic space:
\begin{equation}
    \bar{\bm{F}}_{i\pm\f{1}{2},j} = \int_{K_y^j}  \hat{\bm{F}}_{i\pm\f{1}{2}}  \mu(\bm{y}) \: d\bm{y}.
    \label{eq:Fbar}
\end{equation}

The SFV method then solves the following ODE system:
\begin{equation}
    \f{\rmd \bm{U}_{i,j}} {\rmd\rmt} + \frac{\bar{\bm{F}}_{i+\f{1}{2},j} - \bar{\bm{F}}_{i-\f{1}{2},j}}{\LRv{K_x^i}\LRv{K_y^j}} = \bm{0},
    \label{eq:SFV method_elementwise}
\end{equation}
for all $i=1,...,N_x$, $j=1,...,N_y$. 

We can also write the SFV method in matrix form. Let $\bar{\bm{F}}_{\pm\f{1}{2}}$ denote matrices with entries 
\[
\LRp{\bar{\bm{F}}_{\pm\f{1}{2}}}_{i,j} = \bar{\bm{F}}_{i\pm\f{1}{2},j}.
\]
We can then rewrite \eqref{eq:SFV method_elementwise} as
\begin{equation}
    \bm{M}\f{\rmd \bm{U}} {\rmd\rmt} + \LRp{\bar{\bm{F}}_{+\f{1}{2}} - \bar{\bm{F}}_{-\f{1}{2}}} = \bm{0},
    \label{eq:SFV method_global}
\end{equation}
where $\bm{U}$ is the matrix of cell averages, $\bar{\bm{F}}_{\pm \f{1}{2}}$ denotes the flux integrals over the physical cell boundaries, and $\bm{M}$ is a mass matrix such that
\[
\LRp{\bm{M} \f{\rmd \bm{U}} {\rmd\rmt}}_{i,j} = \LRv{K_x^i}\LRv{K_y^j} \f{\rmd \bm{U}_{i,j}} {\rmd\rmt}.
\]

\begin{remark}
While we focus on one-dimensional spatial domains and higher dimensional stochastic domains in this work, the SFV method framework extends naturally to systems in two/three-dimensional spatial domains \cite{Tokareva14}.
\end{remark}

\section{WENO reconstruction}
\label{sec:WENO}
To achieve high-order accuracy, we use a third-order WENO scheme \cite{Liu94,Titarev04} by applying WENO reconstructions over both the physical and stochastic spaces \cite{Tokareva14,Abgrall2017}. 

\subsection{Formulation and the physical space reconstruction}
We start by presenting the details for a third-order WENO reconstruction. Suppose there exists a cell with cell average $U_i$ connected to its neighbors, we define a three-point stencil with its neighboring solution averages $U_{i\pm1}$. Let $a \in [-1/2,1/2]$ be the reference location within the cell where we want to evaluate the reconstructed value $\Tilde{U}$. We define two candidate linear reconstructions in terms of $a$ 
\begin{equation}
p_0(a) = U_i + a\LRp{U_{i+1}-U_i}, \quad p_1(a) =  U_i + a \LRp{U_i-U_{i-1}}.
    \label{eq:WENO_p}
\end{equation}
Then the reconstructed value at position $a$ is computed as a weighted combination of candidate reconstructions at the cell interfaces
\begin{equation}
    \Tilde{U} = \omega_0 p_0(a) + \omega_1 p_1(a).
    \label{eq:WENO_U}
\end{equation}
To adaptively determine the contribution of these candidate reconstructions and assign higher weights to smoother stencils, we compute the smoothness indicators
\[
\beta_0 = \LRp{U_{i+1}-U_i}^2\csp \beta_1 = \LRp{U_i - U_{i-1}}^2.
\]
The corresponding nonlinear weights are given by
\begin{equation}
    \alpha_k = \f{d_k}{\LRp{\epsilon + \beta_k}^2}\csp \omega_k = \f{\alpha_k}{\alpha_0+\alpha_1}\csp k = 0,1,
    \label{eq:WENO_weights}
\end{equation}
where $\epsilon$ is a small parameter to prevent division by zero. The linear optimal weights $d_0$ and $d_1$ are chosen so that the reconstruction exactly reproduces any quadratic polynomial in the reference cell
\begin{equation}
d_0 = \frac{a^2 + a - 1/12}{2a} \csp
d_1 = 1 - d_0.
\label{eq:WENO_d}
\end{equation}

We propose two different procedures for reconstructing the numerical flux: WENO with reconstructed states or with reconstructed fluxes. Both of these procedures begin by applying third-order WENO reconstruction in the physical domain to obtain the reconstructed left and right solution states $\Tilde{\bm{u}}_{L,i\pm\f{1}{2},j}$, $\Tilde{\bm{u}}_{R,i\pm\f{1}{2},j}$ at the physical boundaries of the cell $K_x^i.$ We first obtain $d_0 = 1/3, d_1=2/3$ for the left interface ($a=-1/2$) $x_{i-\f{1}{2}}$ and $d_0 = 2/3, d_1=1/3$ for the right interface ($a=1/2$) $x_{i+\f{1}{2}}$ from \eqref{eq:WENO_d}. The reconstructed states, $\Tilde{u}_{R,i-\f{1}{2},j}$ and $\Tilde{u}_{L,i+\f{1}{2},j}$, are then computed by plugging $a$ into \eqref{eq:WENO_p} and \eqref{eq:WENO_U}
\begin{equation}
    \begin{split}
         & \Tilde{u}_{R,i-\f{1}{2},j} = \omega_0 \LRp{U_{i,j} - \f{1}{2}\LRp{U_{i+1,j}-U_{i,j}}} + \omega_1 \LRp{U_{i,j}-\f{1}{2}\LRp{U_{i,j}-U_{i-1,j}}}, \\
    &\Tilde{u}_{L,i+\f{1}{2},j} = \omega_0 \LRp{U_{i,j} + \f{1}{2}\LRp{U_{i+1,j}-U_{i,j}}} + \omega_1 \LRp{U_{i,j}+\f{1}{2}\LRp{U_{i,j}-U_{i-1,j}}}.
    \end{split}
    \label{eq:WENO_1D}
\end{equation}
The remaining two reconstructed states, $\Tilde{u}_{R,i+\f{1}{2},j}$ and $\Tilde{u}_{L,i-\f{1}{2},j}$, can be obtained by applying the same procedure to physically neighboring cells. 

\begin{remark}
   The formal third-order accuracy of the WENO scheme is achieved only for smooth solutions; in the presence of shocks, the order of accuracy drops to first order \cite{Tokareva14,Abgrall2017}.
\end{remark}

\begin{remark}
  For $n$-component solutions, WENO reconstruction can be applied independently to each component. In the case of systems in two- or three-dimensional physical space, the reconstruction can be performed dimension-by-dimension, as described in \cite{Titarev04}.
\end{remark}

\subsection{WENO with reconstructed states}
\label{sec:state_recons}
Assuming that we have obtained the reconstructed states $\Tilde{\bm{u}}_{L,i\pm\f{1}{2},j}$, $\Tilde{\bm{u}}_{R,i\pm\f{1}{2},j}$ after \eqref{eq:WENO_1D}, the first procedure, as used in \cite{Tokareva14,Abgrall2017}, is WENO with ``reconstructed states'' that extends the reconstruction process over the stochastic domain. 

This is also performed in a cell-by-cell and dimension-by-dimension manner, following a similar procedure as in \eqref{eq:WENO_1D}. As a result, we obtain interpolation polynomials $\Tilde{\bm{u}}_{L,i\pm\f{1}{2},j}\LRp{\bm{y}}$, $\Tilde{\bm{u}}_{R,i\pm\f{1}{2},j}\LRp{\bm{y}}$ for each $K_x^i\times K_y^j$, which enable the evaluation of approximated solution states at arbitrary quadrature nodes. 

For any chosen quadrature rule within one stochastic cell, let $\bm{\hat{y}}_m$ and $w_{m}$ be the corresponding quadrature nodes and weights, respectively. \footnote{Note that quadrature weights $w_m$ are distinct from WENO weights $\omega_k$.} Using this quadrature rule, we approximate the flux integrals in \eqref{eq:Fbar} as
\begin{equation}
 \int_{K_y^j}  \hat{\bm{F}}_{i\pm\f{1}{2}}\LRp{\bm{y}} \mu(\bm{y}) \: d\bm{y}  \approx \sum_{m}\hat{\bm{F}}\LRp{\Tilde{\bm{u}}_{L,i\pm\f{1}{2},j}\LRp{\bm{\hat{y}}_{m}}, \Tilde{\bm{u}}_{R,i\pm\f{1}{2},j}\LRp{\bm{\hat{y}}_{m}}}\mu(\bm{\hat{y}}_{m})w_{m}.
\label{eq:integral_state}
\end{equation}
Since we employ a third-order WENO reconstruction, integrals of numerical fluxes are approximated using two-point Gauss quadrature rules which are exact for cubic polynomials. 

\subsection{WENO with reconstructed fluxes}
\label{sec:WENO_flux}
The second procedure, WENO with ``reconstructed fluxes'', first computes the flux at each cell interface using the reconstructed states: 
\[
\hat{\bm{F}}_{i\pm\f{1}{2},j}\LRp{\bm{y}} = \hat{\bm{F}}\LRp{\Tilde{\bm{u}}_{L,i\pm\f{1}{2},j}, \Tilde{\bm{u}}_{R,i\pm\f{1}{2},j}}.
\]
Then, rather than reconstructing the solution states over the stochastic domain as in the previous procedure, we directly reconstruct the computed flux, following the approach commonly used in finite difference WENO schemes since the original formulation in \cite{JIANG1996}. This process yields the flux interpolation polynomials $\Tilde{\bm{F}}_{i\pm \f{1}{2},j}\LRp{\bm{y}}$ for each $K_x^i\times K_y^j$, which can be evaluated at quadrature nodes. Using the same quadratures described in Section \ref{sec:state_recons}, the flux integrals in \eqref{eq:Fbar} are then approximated as 
\begin{equation}
 \int_{K_y^j}  \hat{\bm{F}}_{i\pm\f{1}{2}}\LRp{\bm{y}} \mu(\bm{y}) \: d\bm{y}  \approx \sum_{m} \Tilde{\bm{F}}_{i\pm \f{1}{2},j}\LRp{\bm{\hat{y}}_{m}} \mu(\bm{\hat{y}}_{m})w_{m}.
\label{eq:integral_flux}
\end{equation}

A key distinction of this method is that it evaluates the numerical flux only in the physical domain, whereas the reconstructed states procedure requires flux evaluations at every quadrature node in the stochastic domain. Specifically, at each ODE time step, the procedure of reconstructed states requires $2^q$ times the flux evaluations compared to the procedure of reconstructed fluxes. Such reduction comes from the fact that two-point Gaussian quadrature rule yields $2^q$ quadrature points in $q$ dimensions. This implies if the flux evaluation is a significant part of the computational cost compared to the reconstruction, the procedure of reconstructed fluxes will be more efficient. Although the total number of time steps depends on the specific time integrator and CFL condition \cite{CFL}, we observe that this procedure reduces the overall runtime of the ODE time-stepping process in practice.

\section{Reduced stochastic basis}
\label{sec:reduced}
Despite potential runtime improvements with the procedure of reconstructed fluxes, both WENO reconstruction procedures still require obtaining the flux at every quadrature node, resulting in an exponential increase in runtime and memory as the stochastic dimension $q$ grows. To alleviate this, we apply interpolation-based reduced order model (ROM) techniques to reduce costs associated with high dimensional stochastic domains. More specifically, we non-intrusively construct a data-driven reduced basis in the stochastic domain, which allows us to approximate high dimensional stochastic integrals in terms of precomputed operators acting on the reduced stochastic basis.

\subsection{Approximating stochastic integrals}
We approximate stochastic flux integrals in \eqref{eq:Fbar}  by representing the flux using $N$ reduced stochastic basis functions $\bm{\phi}(\bm{y}) = [\phi_1(\bm{y}),...,\phi_N(\bm{y})]^T$: %\ST{Explicit form of these stochastic basis functions? Polynomials?}
\[
\hat{\bm{F}}_{i\pm\f{1}{2}}\LRp{\bm{y}} \approx \sum_{k=1}^N \bm{c}_{i\pm\f{1}{2},k} \phi_k(\bm{y}),
\]
where $\bm{c}_{i\pm\f{1}{2},k}\in \mathbb{R}^n$ is a vector of coefficients for the $k$-th basis at the physical interfaces $i\pm\f{1}{2}$; each entry in this vector corresponds to the same solution component as the respective entry in $\hat{\bm{F}}_{i\pm\f{1}{2}}$. The explicit forms of the basis functions can be polynomials or data-driven basis functions such as POD modes, which will be discussed in Sec \ref{sec:POD}. Substituting this expansion into \eqref{eq:Fbar}, the flux integral transforms into 
\begin{equation}
\begin{split}
    \int_{K_y^j}  \hat{\bm{F}}_{i\pm\f{1}{2}}\LRp{\bm{y}} \mu(\bm{y}) \: d\bm{y}  & \approx  \int_{K_y^j}  \sum_{k=1}^N \bm{c}_{i\pm\f{1}{2},k} \phi_k(\bm{y}) \mu(\bm{y}) \: d\bm{y}  \\
    & =  \sum_{k=1}^N \bm{c}_{i\pm\f{1}{2},k} \int_{K_y^j}  \phi_k(\bm{y}) \mu(\bm{y}) \: d\bm{y}.
\end{split}
    \label{eq:integral_interp}
\end{equation}

We use the same notations for quadrature rule as in Section \ref{sec:state_recons}.  Let $\bm{B}\in\mathbb{R}^{N_y\times N}$ be the matrix representing the integrals over stochastic faces of the basis functions weighted by the density functions:
\begin{equation}
\bm{B}_{j,k} = \int_{K_y^j}  \phi_k(\bm{y}) \mu(\bm{y}) \: d\bm{y} \approx \sum_{m}\phi_k(\bm{\hat{y}}_{m}) \mu(\bm{\hat{y}}_{m})w_{m}.
\label{eq:Bmatrix}
\end{equation}
This allows us to trade the online cost of a high-dimensional stochastic integral for offline computations of integrals of reduced stochastic basis functions.

For each spatial interface between neighboring cells, denote $N_q$ as the number of quadrature nodes per stochastic cell, the total number of quadrature nodes along such spatial interface becomes $N_q N_y$, which is typically large. We use a least-squares approximation to determine the coefficients $\bm{c}_{i\pm\f{1}{2},k}$ of the reduced stochastic basis flux approximation. Instead of using $\bm{\hat{y}}_m$ to refer to a local quadrature node within a stochastic cell, we now introduce $\bm{y}_l$ to represent the global $l$-th quadrature node, where $l=1,...,N_qN_y.$

Recall from \eqref{eq:integral_flux} that $\Tilde{\bm{F}}_{i\pm \f{1}{2},j}\LRp{\bm{y}}$ represent the reconstructed flux polynomials that can be evaluated at each quadrature node in $K_x^i\times K_y^j$. To organize these values, we define matrices $\Tilde{\bm{F}}_{\pm\f{1}{2}}$ which store the reconstructed flux at all quadrature nodes in $D_x\times D_y$:
\begin{equation}
    \LRp{\Tilde{\bm{F}}_{\pm\f{1}{2}}}_{l,i} = \Tilde{\bm{F}}_{i\pm \f{1}{2},j}\LRp{\bm{y}_{l}}\csp l = 1,...,N_qN_y,\: i = 1,...,N_x,
    \label{eq:F_global}
\end{equation}
where $j$ corresponds to the cell index such that $N_y^j$ contains the quadrature node $l$. We define the corresponding Vandermonde matrix $\bm{V} \in \mathbb{R}^{N_qN_y\times N}$ as
\begin{equation}
    \bm{V}_{l,k} = \phi_k(\bm{y}_{l}).
    \label{eq:V}
\end{equation}
To approximate the values of $\bm{c}_{i\pm\f{1}{2},k}$, we solve a standard linear least square problem:
\begin{equation}
\min_{\bm{P} \in \text{span}\{\bm{\phi}_1, \dots, \bm{\phi}_N\}} \sum_l \left\| \Tilde{\bm{F}}_{i\pm \frac{1}{2},j}(\bm{y}_l) - \bm{P}(\bm{y}_l) \right\|_2^2.
\label{eq:lsq}
\end{equation}
The analytical solution of the problem above is:
\[
\bm{P} = \bm{V}\bm{V}^\dagger\Tilde{\bm{F}}_{\pm\f{1}{2}},
\]
where $\dagger$ denotes the pseudoinverse. In particular, note that $\bm{V}^\dagger\Tilde{\bm{F}}_{\pm\f{1}{2}}$ are the coefficients of the solution of \eqref{eq:lsq} in the basis $\bm{\phi}_k(\bm{y})$. Since $\bm{B}$ in \eqref{eq:Bmatrix} maps coefficients in the $\bm{\phi}_k(\bm{y})$ basis to integrals over stochastic faces, $\bar{\bm{F}}_{\pm\f{1}{2}}$ in \eqref{eq:SFV method_global} can be approximated as
\begin{equation}
    \bar{\bm{F}}_{\pm\f{1}{2}} \approx \LRp{\bm{B}\bm{V}^\dagger\Tilde{\bm{F}}_{\pm\f{1}{2}}}^T.
\end{equation}
This formulation enables us to approximate the solutions in \eqref{eq:SFV method_global} by solving the reduced-order interpolation-based model
\begin{equation}
    \bm{M}\f{d \bm{U}} {dt} + \LRp{\bm{B}\bm{V}^\dagger\Tilde{\bm{F}}_{+\f{1}{2}} - \bm{B}\bm{V}^\dagger\Tilde{\bm{F}}_{-\f{1}{2}}}^T = \bm{0}.
    \label{eq:ROM}
\end{equation}
Note that even with a reduced basis approach, the dimensions of $\bm{V}$ and $\Tilde{\bm{F}}_{\pm\f{1}{2}}$ still scale with the size of the stochastic space. As a result, the computational cost of \eqref{eq:ROM} does not necessarily decrease compared to the full SFV method.

\subsection{Q-DEIM hyper-reduction}
To scale the dimensions of $\bm{V}$ and $\Tilde{\bm{F}}_{\pm\f{1}{2}}$ with the number of reduced basis, a further step of hyper-reduction on \eqref{eq:ROM} is performed to select a reduced set of interpolation points from Q-DEIM \cite{Drmac16}.

Specifically, given $\bm{V}$ in \eqref{eq:V}, we obtain a permutation vector $\bm{\pi}$ from the QR factorization of $\bm{V}^T.$ Subsequently, we pick $N_H$ leading entries of $\bm{\pi}$ and denote them as
\[
\mathcal{I} = \bm{\pi}[1:N_H].
\]
In this way, $\mathcal{I}$ picks the leading rows of $\bm{V}$ that are most significant with largest norm. We then obtain a hyper-reduced version of \eqref{eq:ROM}:
\begin{equation}
     \bm{M}\f{d \bm{U}} {dt} + \LRp{\bm{B}\bm{V}\LRs{\mathcal{I},:}^\dagger\Tilde{\bm{F}}_{+\f{1}{2}}\LRs{\mathcal{I},:} - \bm{B}\bm{V}\LRs{\mathcal{I},:}^\dagger\Tilde{\bm{F}}_{-\f{1}{2}}\LRs{\mathcal{I},:}}^T = \bm{0}.
    \label{eq:HR_ROM}
\end{equation}

To summarize, in the full SFV method, a total number of $N_q N_y$ quadrature points are used on each physical interface to evaluate the flux integrals over the stochastic cells. In the ROM, we interpolate the fluxes over all available quadrature points with $N$ reduced basis functions and approximate each flux integral over the stochastic cells with some linear combination of the integrals of these $N$ interpolation bases. In the hyper-reduced ROM, we use Q-DEIM to select a subset of $N_H$ interpolation points from the full set in ROM. In each time-step, all flux integrals can therefore be efficiently approximated by evaluating the nonlinear fluxes only at these $N_H$ selected interpolation points to compute the corresponding interpolation coefficients.

To evaluate the reconstructed fluxes at the hyper-reduced indices \(\mathcal{I}\), one could naively compute the full matrices \(\Tilde{\bm{F}}_{\pm\f{1}{2}}\) and then extract the rows  corresponding to $\mathcal{I}$. However, this approach eliminates any computational advantage, making the ROM no more efficient than the full SFV method. In practice, for each quadrature point indexed by $\mathcal{I}$, we instead construct a minimal local stencil consisting of its neighboring stochastic cells in each stochastic dimension. This leads to a reduced tensor-product stencil in the stochastic space. 

\begin{figure}
\centering
  \includegraphics[width=0.8\textwidth]{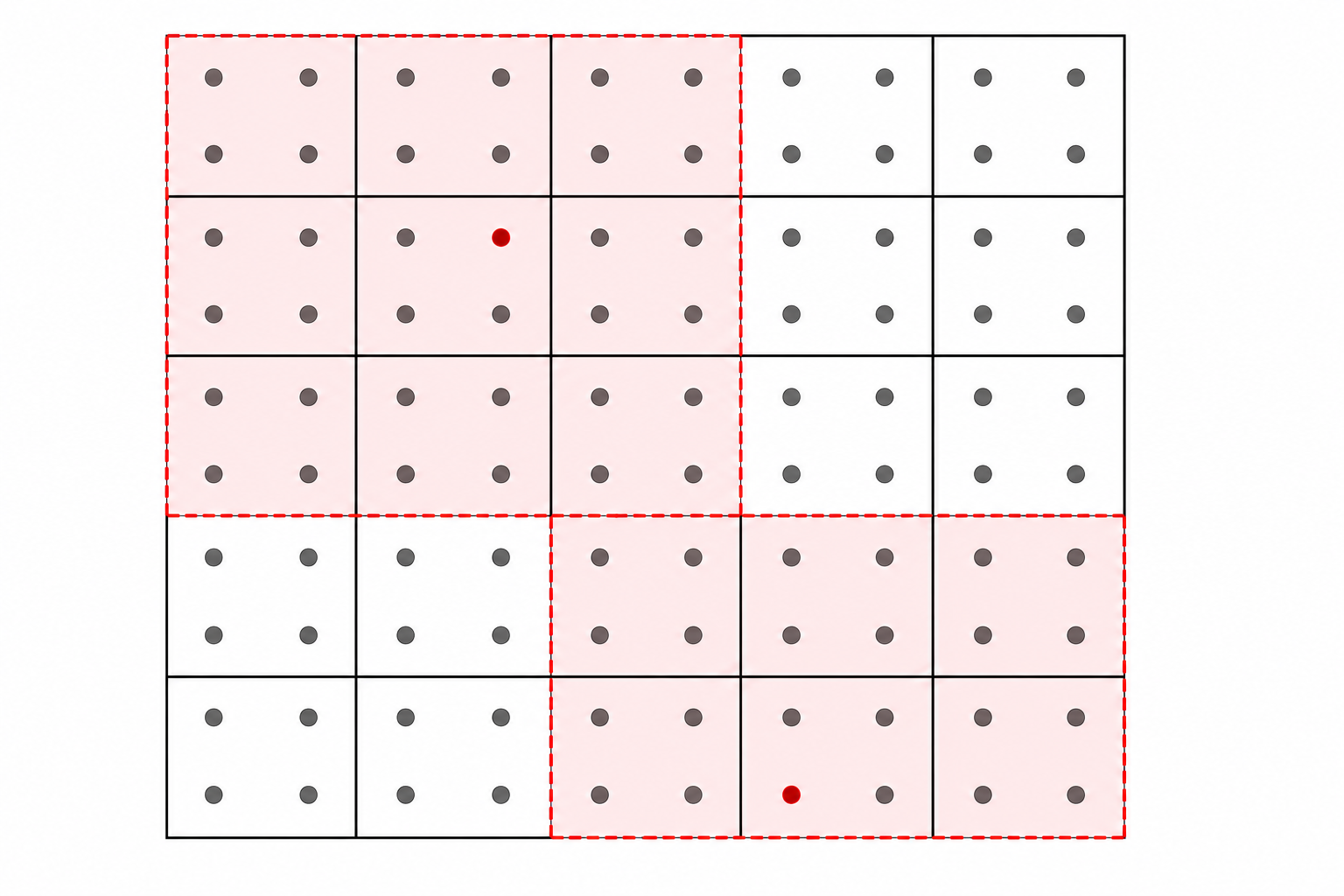}
\caption{Illustration of the hyper-reduction procedure on a $5\times 5$ stochastic grid with two selected quadrature/interpolation points and their associated local stencils highlighted in red.}
\label{fig:HR_stencil}      
\end{figure}

To better illustrate this process, we include a schematic in \autoref{fig:HR_stencil}, where we consider a $5\times 5$ discretization of a two-dimensional stochastic space, employ two-point Gaussian quadrature, and select two hyper-reduced interpolation points. In this setting, it suffices to evaluate the physical fluxes only on the highlighted local stencils prior to applying WENO reconstruction to get the reconstructed fluxes at the hyper-reduced points. 

To further improve efficiency, WENO reconstruction in each stochastic dimension is performed only on the subset of indices associated with the selected quadrature points in that coordinate, rather than over the entire reduced stochastic grid. This avoids redundant computations of fluxes that are not used. Nevertheless, we note that the total number of nonlinear WENO reconstructions still scales exponentially with the stochastic dimension $q$.

\subsection{Choice of reduced stochastic basis}
\label{sec:POD}
Our hyper-reduced model should work with any choice of reduced stochastic basis as long as it is not ill-conditioned (for example, we could use a polynomial basis in stochastic space). % Legendre polynomials \cite{Mccarthy93} and Chebyshev polynomials \cite{Mason02}. 
For this work, we utilize a non-intrusively constructed reduced stochastic basis using POD \cite{Chatterjee00,Liang01}.

Recall that $\Tilde{\bm{F}}_{\pm\f{1}{2}}$ store the reconstructed flux at all quadrature nodes, as defined in \eqref{eq:F_global}. Now, let $\Tilde{\bm{F}}_{\pm\f{1}{2},p} \in \mathbb{R}^{N_qN_y\times N_x}$ represent the reconstructed flux of the $p$th solution component of the conservation laws. Then, for a set of equally spaced time frames $[t_1,...,t_m]$, we construct
\begin{equation}
\label{eq:Fsnap}
\begin{array}{c}
\LRp{\Tilde{\bm{F}}_\text{snap}}_{\pm\f{1}{2}} = \begin{bmatrix}
    \Tilde{\bm{F}}_{\pm\f{1}{2},1}(t_1) & \hdots & \Tilde{\bm{F}}_{\pm\f{1}{2},n} (t_1) & \hdots & \Tilde{\bm{F}}_{\pm\f{1}{2},1}(t_m) & \hdots & \Tilde{\bm{F}}_{\pm\f{1}{2},n}(t_m)  
\end{bmatrix}, \\ 
\Tilde{\bm{F}}_\text{snap} = \begin{bmatrix}
    \LRp{\Tilde{\bm{F}}_\text{snap}}_{+\f{1}{2}} & \LRp{\Tilde{\bm{F}}_\text{snap}}_{-\f{1}{2}}
\end{bmatrix} \in \mathbb{R}^{N_qN_y\times 2mnN_x}.
\end{array}
\end{equation}
The reduced basis can be obtained by computing a singular value decomposition (SVD):
\[
\Tilde{\bm{F}}_\text{snap} = \bm{U}_{S}\bm{\Sigma}\bm{V}_{S},
\]
where we choose first $N$ left singular vector modes as our reduced basis. The Vandermonde matrix becomes 
\[
\bm{V} = \bm{U}_S[:,1:N].
\]
We observe that $\Tilde{\bm{F}}_{i\pm \f{1}{2},j}\LRp{\bm{y}}$ are identical to $\Tilde{\bm{F}}_{(i\pm1)\mp \f{1}{2},j}\LRp{\bm{y}}$ except at the boundaries of $D_x$. Consequently, redundant flux snapshots can be omitted to reduce memory usage. Additionally, the randomized SVD \cite{Drinea01,Halko11} can be employed for further memory efficiency. Instead of performing a full SVD on the snapshot matrix, randomized SVD first projects the data onto a lower-dimensional subspace using a random sampling matrix, capturing the dominant column space with high probability. The SVD is then performed on this much smaller projected matrix, significantly reducing both memory usage and computational cost.

All results presented on this work use identical stochastic basis functions shared across all physical interfaces for efficiency. However, using locally adapted stochastic bases could further improve accuracy near shocks, as demonstrated in some adaptive ROM frameworks \cite{AmsallenFarhat}.

\subsection{Non-intrusive construction of reduced stochastic basis}
A key aspect of the approach presented in this paper is that the reduced stochastic basis can be constructed without needing to directly evaluate the full high-dimensional SFV model. In other words, the reduced stochastic basis can be constructed purely by sampling the one-dimensional spatial model at various parameter values. The process of non-intrusively constructing the reduced stochastic basis is straightforward if the reduced stochastic basis is deterministic (e.g., a high order polynomial basis), but can also be done for a data-driven basis such as POD. We outline the process below.

To collect the flux snapshots used in constructing the POD basis, we may solve the full system \eqref{eq:SFV method_global} using any of the reconstruction methods described above. Since the reconstructed fluxes at global quadrature points are explicitly computed at each time step, they can be directly assembled into the snapshot matrix $ \Tilde{\bm{F}}_\text{snap} $ in \eqref{eq:Fsnap}. Alternatively, one can store the solution states $\bm{U}(t_1), ..., \bm{U}(t_m)$ from the full SFV method and later apply reconstruction to these snapshots to recover the corresponding fluxes, which are then used to form $ \Tilde{\bm{F}}_\text{snap}$.

However, the snapshot collection methods above result in an intrusive construction that requires the computation of a full high-dimensional solution. In contrast, we propose a non-intrusive construction by constructing the reduced stochastic basis through a set of decoupled, one-dimensional simulations — one at each stochastic cell center $\bm{y}_j$ of $K_y^j$ for all $j=1,...,N_y$. Specifically, we use the standard finite volume method to solve:
\[
\begin{array}{c}
\partial_t\bm{u} + \nabla_x\cdot\bm{F}\LRp{\bm{u},\bm{y}_j}  = \bm{0}\csp x\in D_x\subset\mathbb{R}\csp t>0;\\
\bm{u}\LRp{x,0,\bm{y}_j} = \bm{u}_0\LRp{x,\bm{y}_j}\csp x\in D_x\in\mathbb{R}.
\end{array}
\]
We first evaluate the flux at the spatial interfaces $i \pm \tfrac{1}{2}$ using reconstruction of the solution snapshots. The WENO reconstruction is then applied to these flux evaluations to obtain the flux values at all quadrature points $\bm{y}_l$, which corresponds to the WENO with reconstructed fluxes scheme in Section~\ref{sec:WENO_flux} and produces the flux snapshot matrix $ \Tilde{\bm{F}}_\text{snap}$ in \eqref{eq:Fsnap}. In this way, the total number of one-dimensional sampling simulations depends only on the Cartesian discretization of the stochastic space, while applying WENO reconstruction to the snapshots significantly enhances accuracy and robustness.

\begin{remark}
    Alternatively, one may sample one-dimensional simulations at each quadrature point $\bm{y}_l$ for all $l = 1,\dots, N_q N_y$ and assemble them to construct $\Tilde{\bm{F}}_{\text{snap}}$. However, this approach requires exponentially more one-dimensional simulations compared to the method we use. Moreover, we observe that the resulting POD modes derived from such $\Tilde{\bm{F}}_{\text{snap}}$ fail to accurately represent the flux evaluation space, leading to large ROM errors or even instabilities.
\end{remark}

\begin{remark}
    One might notice our sampling strategy can already provide complete UQ estimates, as it effectively resolves the stochastic dependence by collecting flux snapshots at every quadrature point. However, one of the main advantages of this framework is its flexibility in choosing how the stochastic space is sampled when constructing the reduced stochastic basis. Instead of sampling for fluxes at all quadrature points, one could select a subset based on sparse grid rules \cite{MA2009,XiuHesthaven05} or other adaptive strategies \cite{Thanh08,Narayan14}, potentially reducing the number of one-dimensional simulations while maintaining sufficient accuracy for UQ purposes.
\end{remark}

\section{Numerical experiments}
\label{sec:num_exp}
In this section, we present numerical experiments to assess the performance of our proposed method for uncertainty quantification across different stochastic dimensions. We consider the one-dimensional compressible Euler equations on some classical problems under different uncertainty settings. We first compare the performance of the two proposed WENO reconstruction procedures and then evaluate the ROM performance.

%order model (ROM).
\subsection{Preliminaries}
We now test the one-dimensional compressible Euler equations:
\begin{equation}
\begin{split}
             & \f{\p \bm{U}}{\p t} + \f{\p \bm{F}(\bm{U})}{\p x} = 0\csp x\in D_x, \\
             & \bm{U} = \LRs{\rho, \rho u, E}^T\csp \bm{F} = \LRs{\rho u, \rho u^2+p, u\LRp{E+p}}^T, \\
             & p = (\gamma - 1)\LRp{E-\f{1}{2}\rho u^2}.
\end{split}
     \label{eq:Euler}
\end{equation}

We choose $\hat{\bm{F}}$ to be the Lax-Friedrich flux \cite{Lax54} which, for  reconstructed left and right states $\Tilde{\bm{u}}_L$, $\Tilde{\bm{u}}_R$, is defined as
\[
\hat{\bm{F}} \LRp{\Tilde{\bm{u}}_L,\Tilde{\bm{u}}_R} = \f{1}{2} \LRp{\bm{F}\LRp{\Tilde{\bm{u}}_L} + \bm{F}\LRp{\Tilde{\bm{u}}_R}} - \f{\lambda}{2}\LRp{\Tilde{\bm{u}}_R - \Tilde{\bm{u}}_L},
\]
where $\lambda$ represents the local maximum wave speed following the Davis wave speed formulation \cite{Davis88}.

On a one-dimensional reference domain $[-1,1]$, the two-point Gaussian quadrature rule approximates integrals via
\[ 
\int_{-1}^1 f(\bm{u})\:d\bm{y} \approx w_1f(y_1) + w_2 f(y_2),
\]
where the quadrature nodes and weights are given by
\[
y_1 = -1/\sqrt{3}\csp y_2 = 1/\sqrt{3}\csp w_1=w_2=1.
\]
For higher-dimensional stochastic integrals, the quadrature rule is constructed as a tensor-product extension of the one-dimensional case. The number of quadrature points in this rule is $N_q = 2^q$, where q denotes the dimension of the stochastic space.

To assess errors of our ROMs, we report the relative $L^1$ norm of the difference between the expected value of the simulated stochastic ROM solutions and that of the reference solution obtained from full SFV method via WENO with reconstructed states or with reconstructed fluxes. The reference solutions are computed on a refined grid consisting of $N_x$ physical volumes, where the fluxes at each interface are evaluated using a total of $N_q N_y$ quadrature points in the stochastic space. Denote $\bm{U} \subset \mathbb{R}^{N_x \times N_q N_y \times n}$ as the computed solutions (cell averages) and denote $\bm{U}_\text{ref}\subset \mathbb{R}^{N_x \times N_qN_y \times n}$ as the corresponding reference solution. Using two-point Gaussian quadrature rule and uniform grid, the relative $L^1$ error simplifies to 
\[
||\bm{U} - \bm{U}_\text{ref}||_{L^1}^{\text{rel}} =
\frac{\sum_{i,j,p} \,\mu(\bm{y}_j)\left|(\bm{U} - \bm{U}_\text{ref})_{i,j,p}\right|}
{\sum_{i,j,p} \,\mu(\bm{y}_j)\left|(\bm{U}_\text{ref})_{i,j,p}\right|}.
\]
where $p$ indexes the solution components.

For all simulations, we use MATLAB's adaptive ODE45 solver with a relative error tolerance of 1e-6 and an absolute error tolerance of 1e-8. The reported runtime is measured using the \texttt{timeit} function applied to the solver and normalized by the total number of time steps. The MATLAB version used is R2024a.

\begin{table}
\caption{Stochastic Sod shock tube 1: relative difference between solutions from full SFV method using WENO with reconstructed states and WENO with reconstructed fluxes, together with their runtime per time step in milliseconds.}
\label{tb:Sod1}
\begin{tabular}{lccc}
\hline\noalign{\smallskip}
$N_y$ & Difference & Runtime per step (state) & Runtime per step (flux) \\
\noalign{\smallskip}\hline\noalign{\smallskip}
$4$  & $5.94\times 10^{-4}$ & $8.12\times 10^{-1}$  & $7.77\times 10^{-1}$ \\
$8$  & $2.87\times 10^{-4}$ & $1.09$  & $1.04$ \\
$16$ & $9.98\times 10^{-5}$ & $1.60$  & $1.58$ \\
$32$ & $5.10\times 10^{-16}$ & $3.09$  & $2.70$ \\
$64$ & $5.39\times 10^{-16}$ & $5.66$ & $5.26$ \\
\noalign{\smallskip}\hline
\end{tabular}
\end{table}

\begin{figure}
        \includegraphics[width=0.5\textwidth]{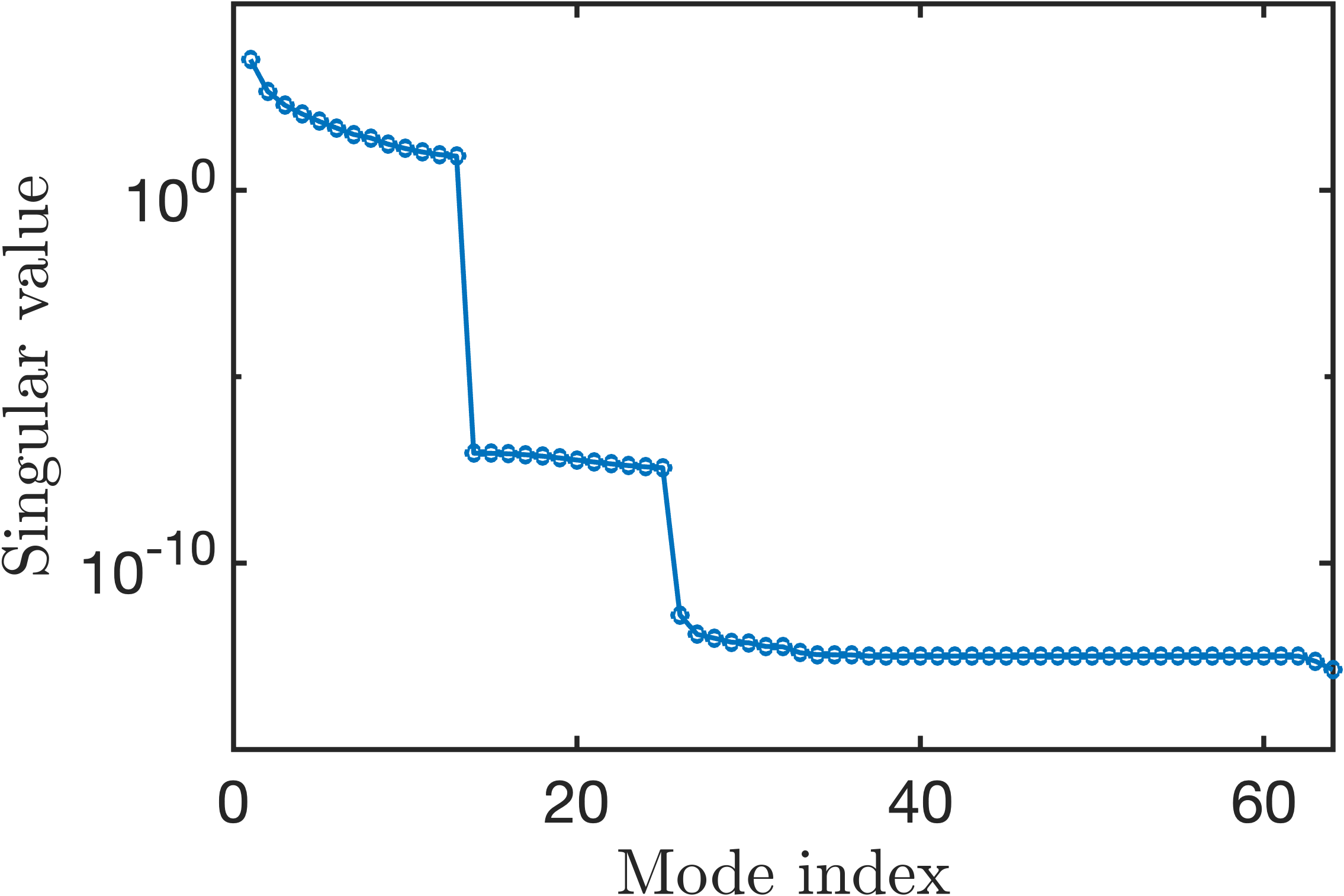}
    \caption{Stochastic Sod shock tube 1: singular values of the flux snapshot matrix.}
    \label{fig:Sod1_svd}
\end{figure}

\begin{figure}
    \centering
    \begin{subfigure}{0.49\textwidth}
        \centering
        \includegraphics[width=\textwidth]{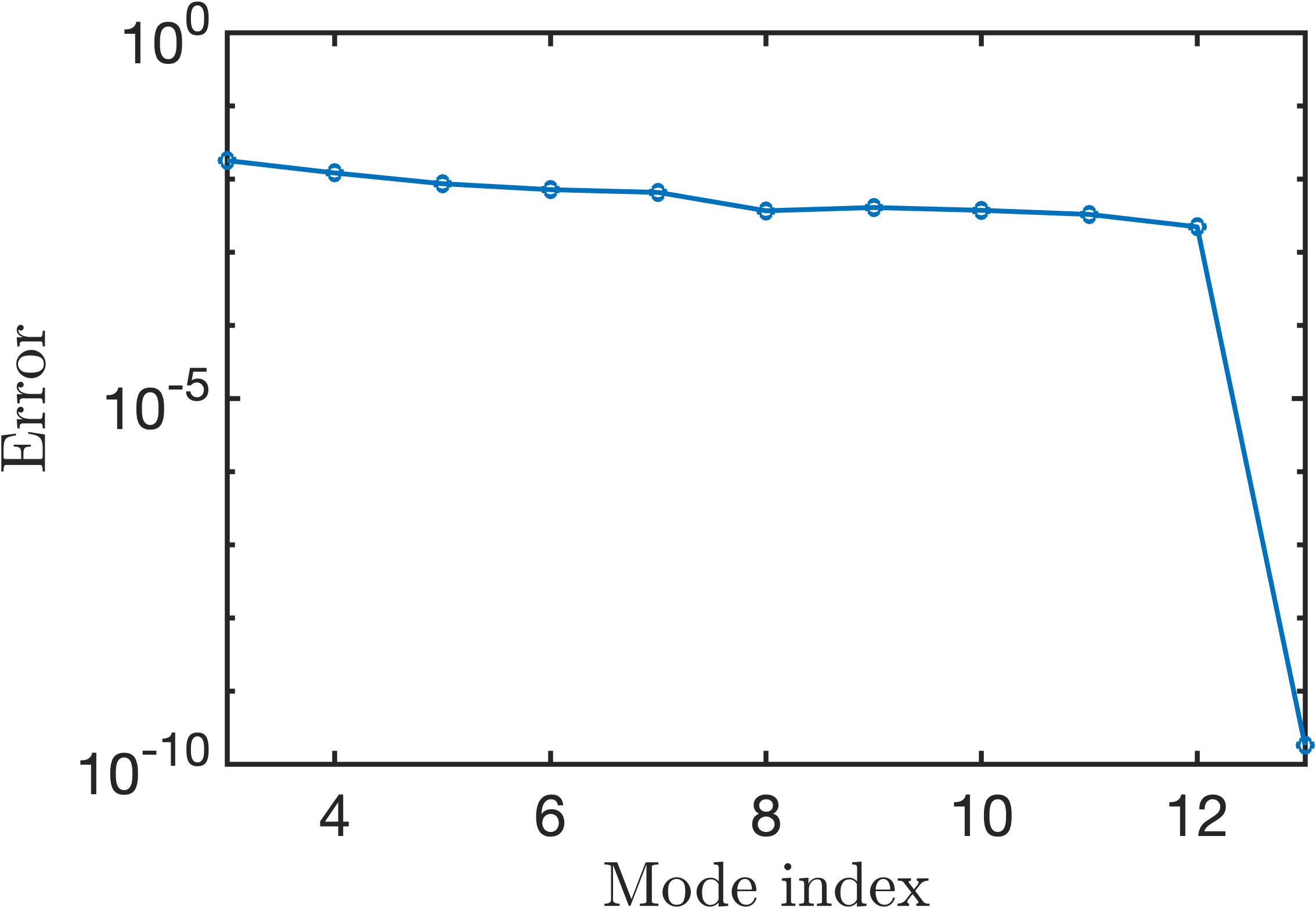}
        \caption{ROM error plot}
    \end{subfigure}
    \begin{subfigure}{0.49\textwidth}
        \centering
        \includegraphics[width=\textwidth]{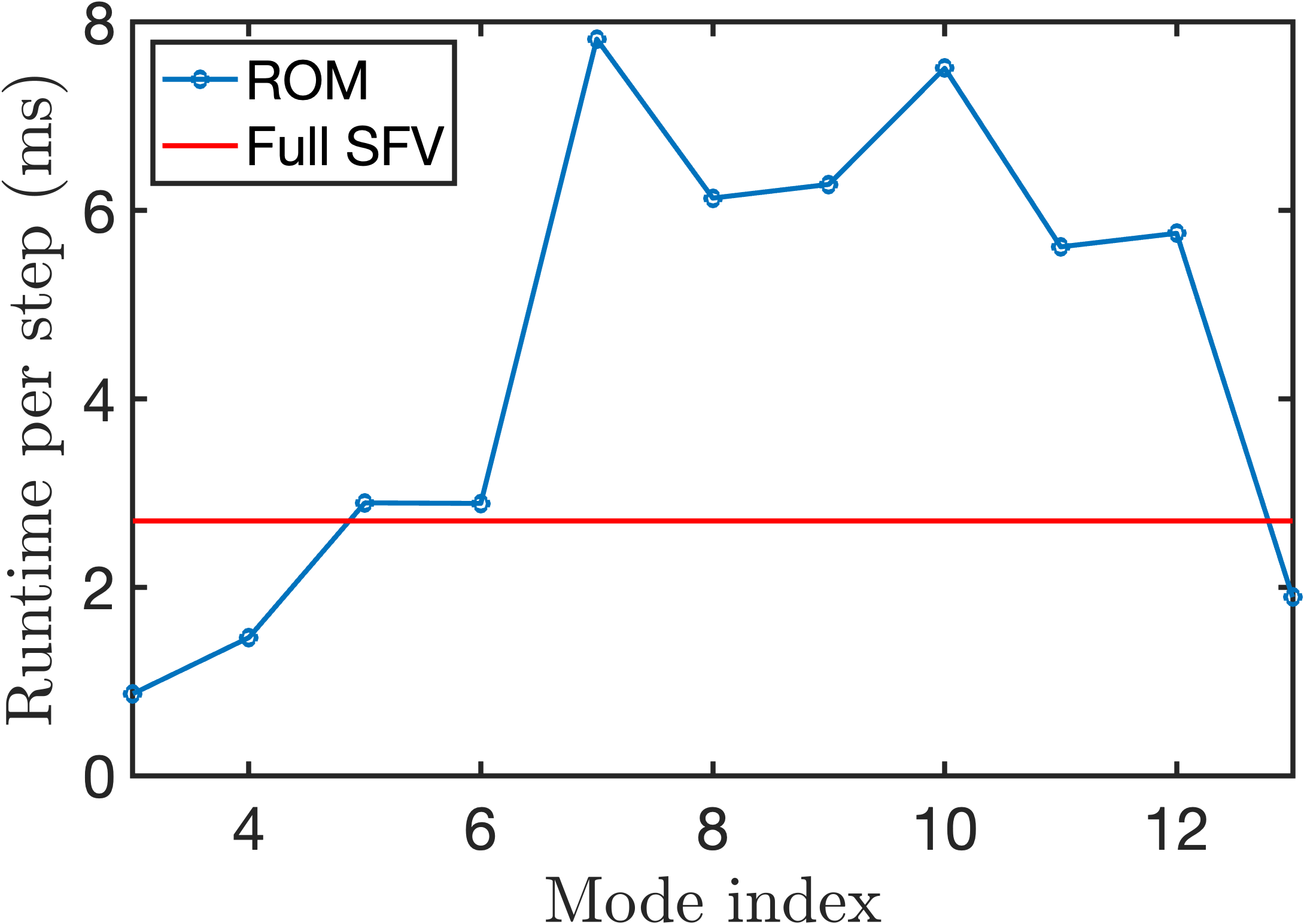}
        \caption{ROM time-stepping runtime plot}
    \end{subfigure}
    \caption{Stochastic Sod shock tube 1: error and runtime per time step for $N_H = N$ hyper-reduced ROMs.}
    \label{fig:Sod1_er_time}
\end{figure}

\begin{figure}
    \centering
    \begin{subfigure}{0.49\textwidth}
        \centering
        \includegraphics[width=\textwidth]{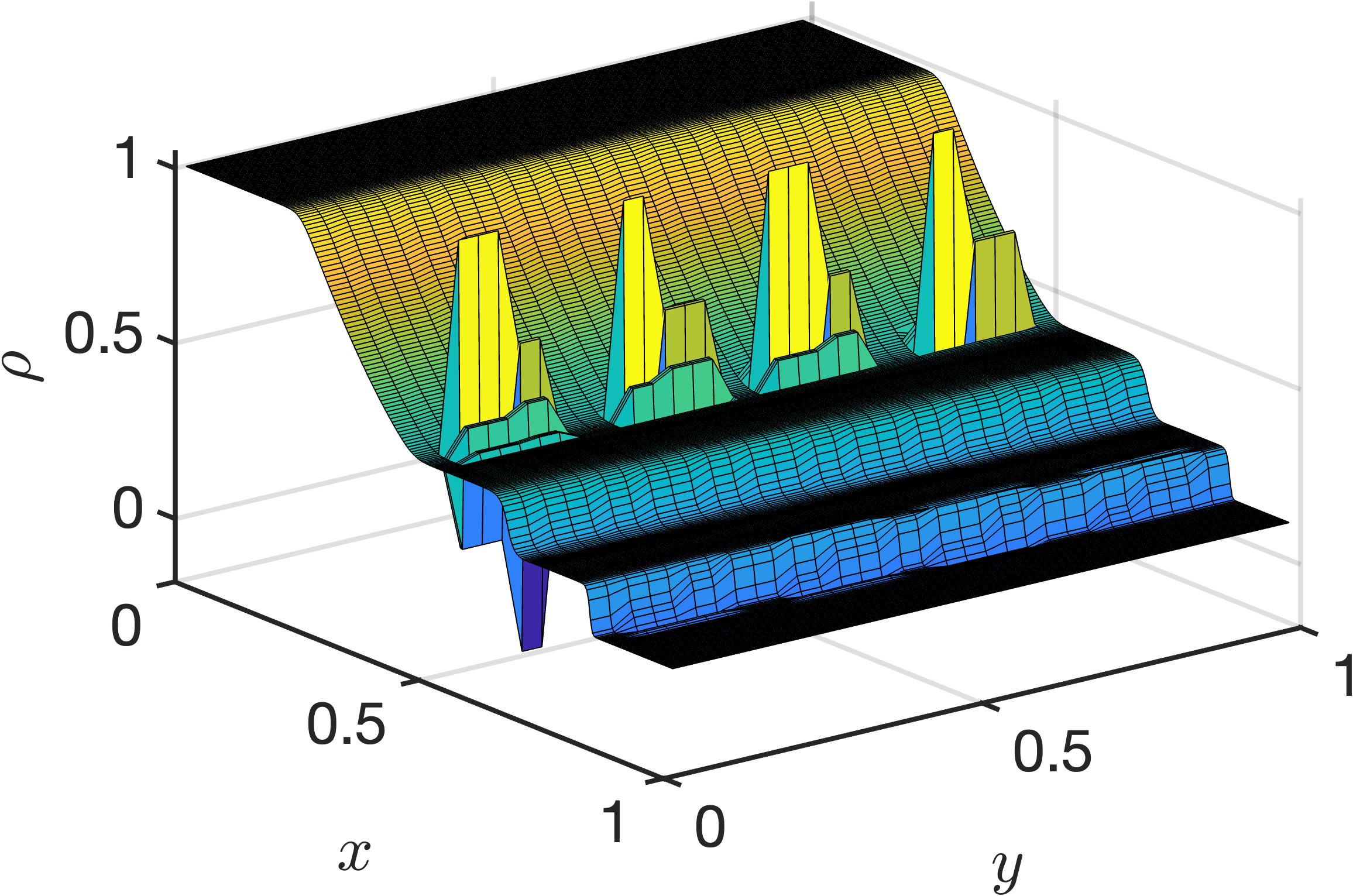}
        \caption{$\rho$ plot, $N=5$ ROM}
    \end{subfigure}
    \begin{subfigure}{0.49\textwidth}
        \centering
        \includegraphics[width=\textwidth]{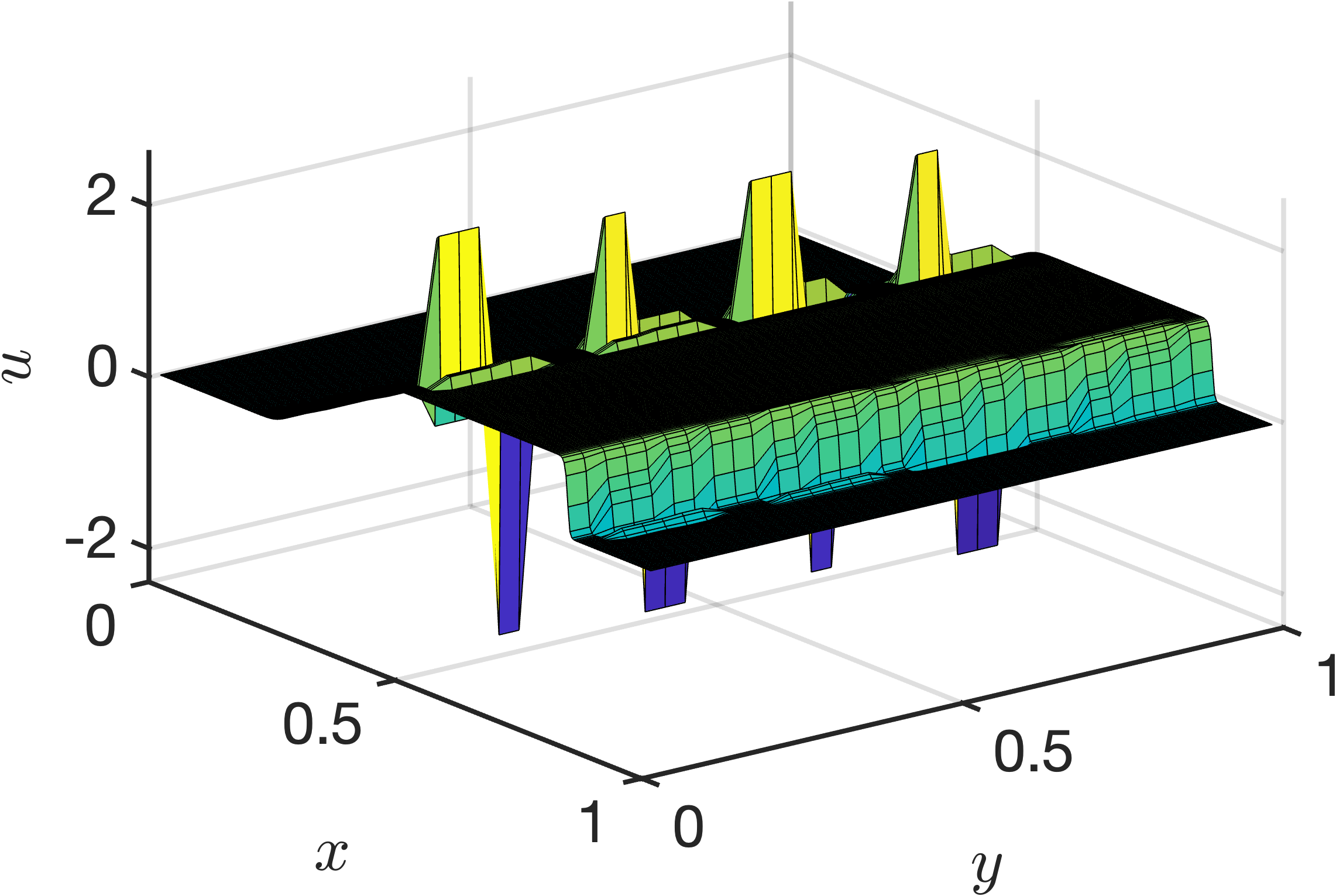}
        \caption{$u$ plot, $N=5$ ROM}
    \end{subfigure}
     \begin{subfigure}{0.49\textwidth}
        \centering
        \includegraphics[width=\textwidth]{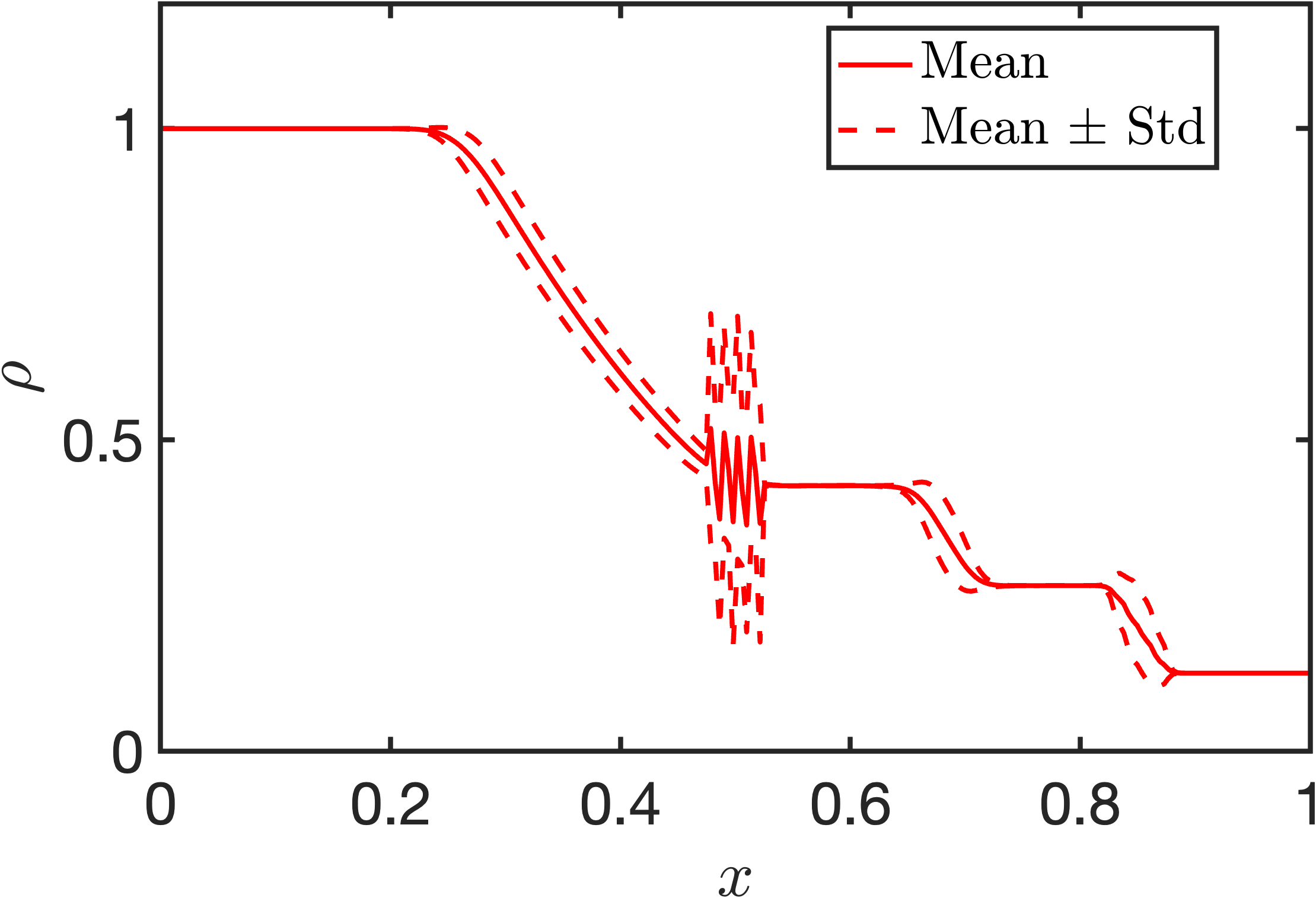}
        \caption{Mean$\pm$std of $\rho$, $N=5$ ROM}
    \end{subfigure}
    \begin{subfigure}{0.49\textwidth}
        \centering
        \includegraphics[width=\textwidth]{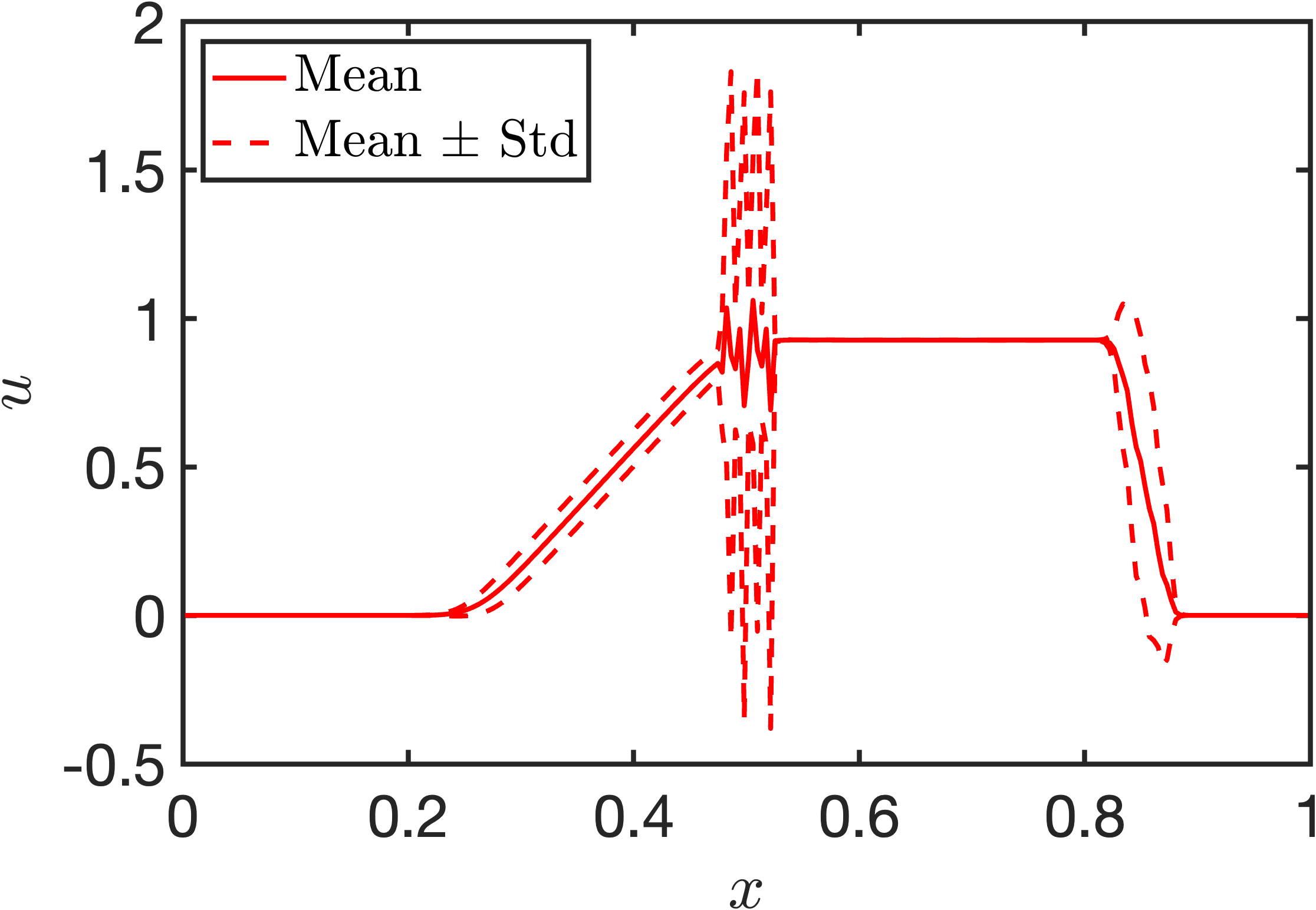}
        \caption{Mean$\pm$std of $u$, $N=5$ ROM}
    \end{subfigure}
    \caption{Stochastic Sod shock tube 1: solution plots and mean for $N_H = N = 5$ (under-resolved) hyper-reduced ROM.}
    \label{fig:Sod1_N5}
\end{figure}

\begin{figure}
    \centering
    \begin{subfigure}{0.49\textwidth}
        \centering
        \includegraphics[width=\textwidth]{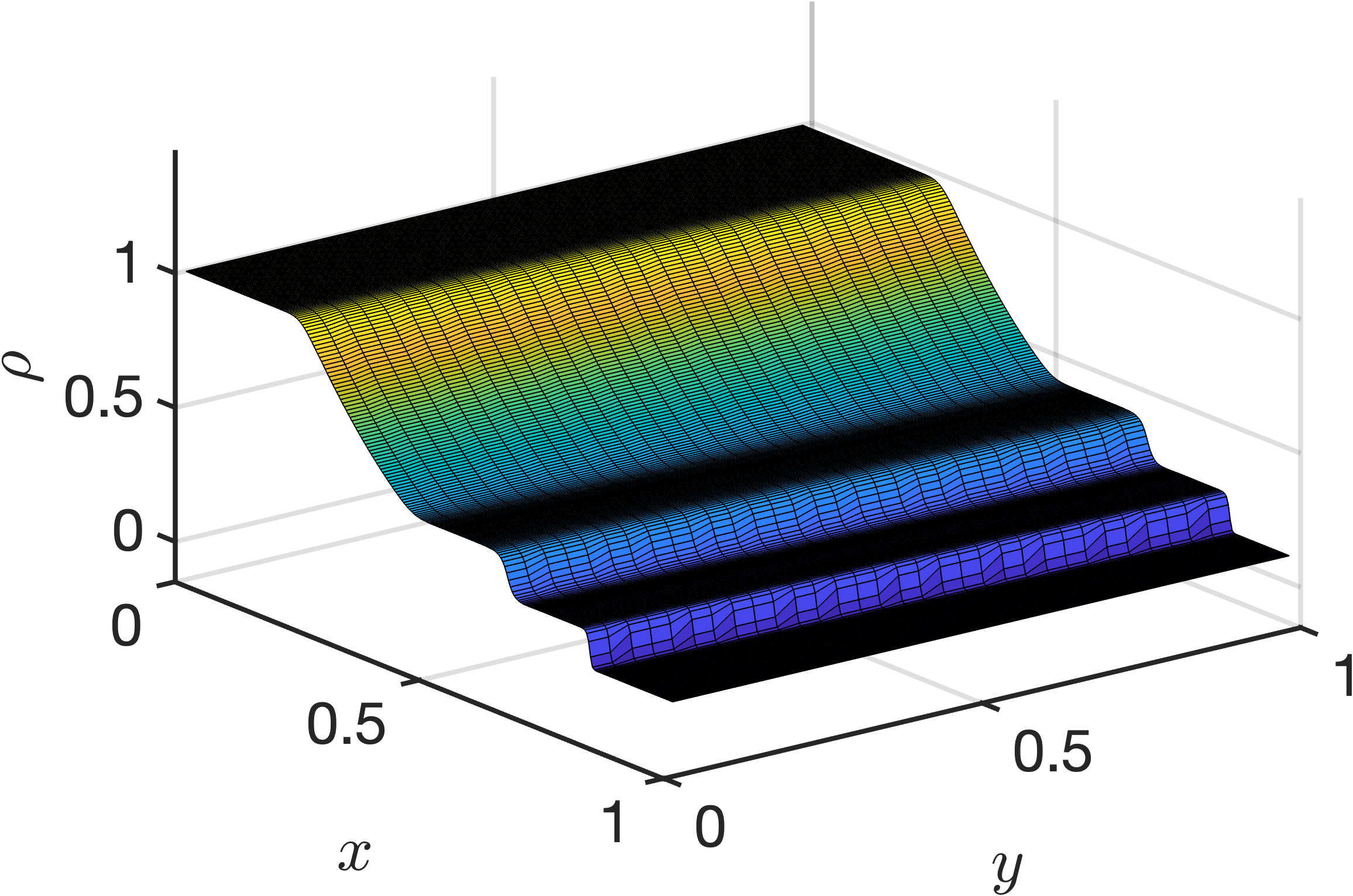}
        \caption{$\rho$ plot, $N=13$ ROM}
    \end{subfigure}
    \begin{subfigure}{0.49\textwidth}
        \centering
        \includegraphics[width=\textwidth]{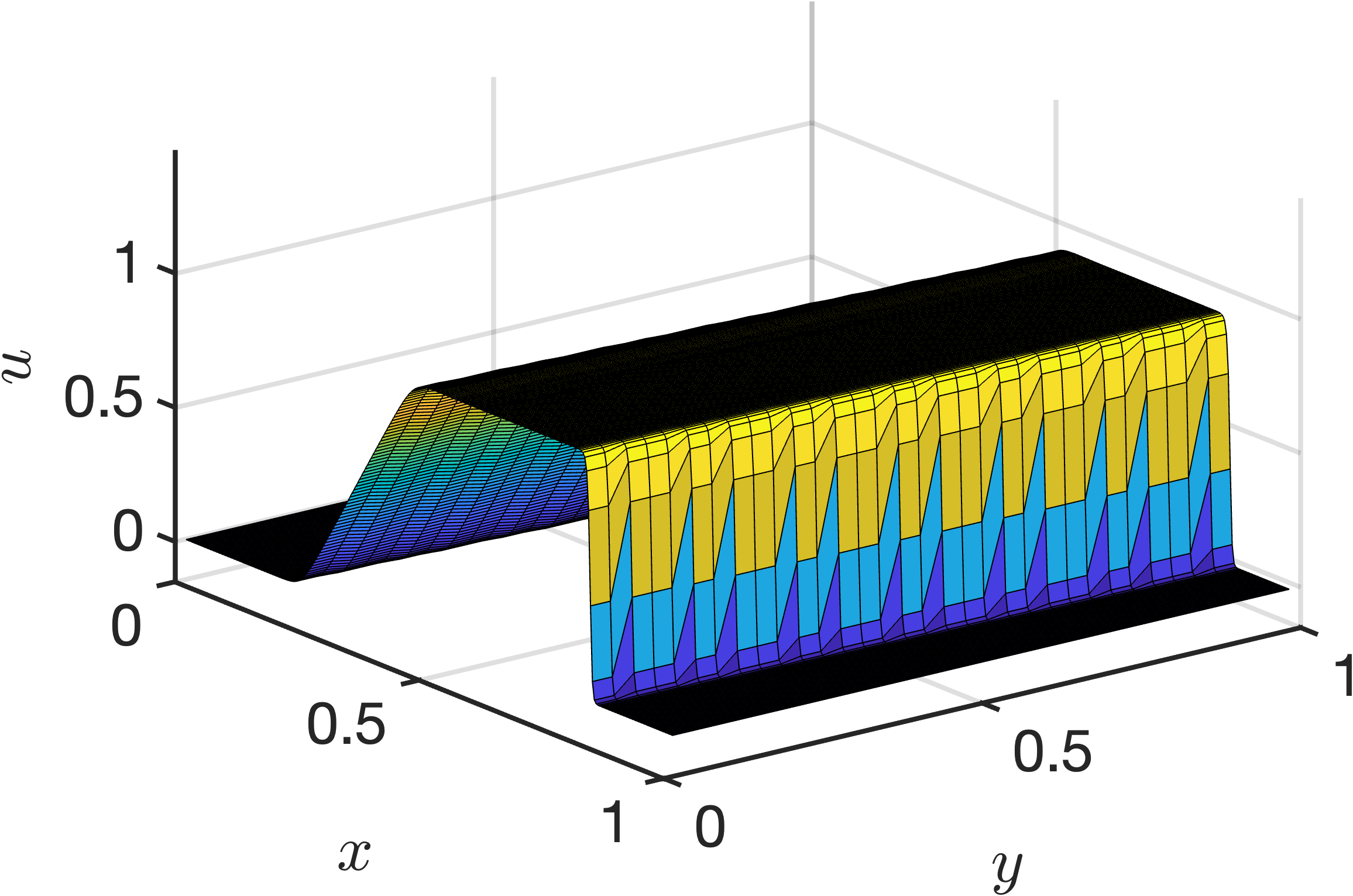}
        \caption{$u$ plot, $N=13$ ROM}
    \end{subfigure}
     \begin{subfigure}{0.49\textwidth}
        \centering
        \includegraphics[width=\textwidth]{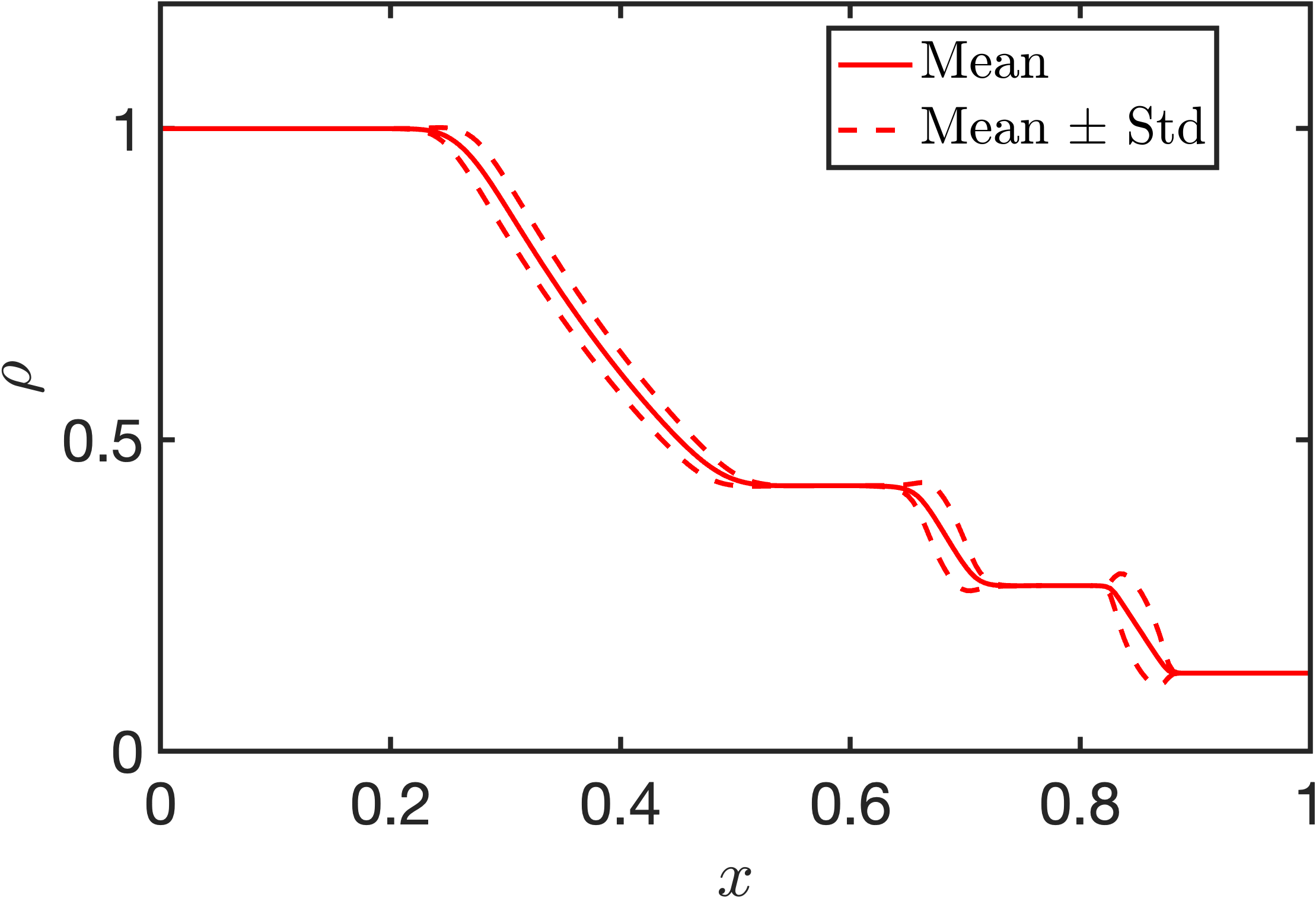}
        \caption{Mean$\pm$std of $\rho$, $N=13$ ROM}
    \end{subfigure}
    \begin{subfigure}{0.49\textwidth}
        \centering
        \includegraphics[width=\textwidth]{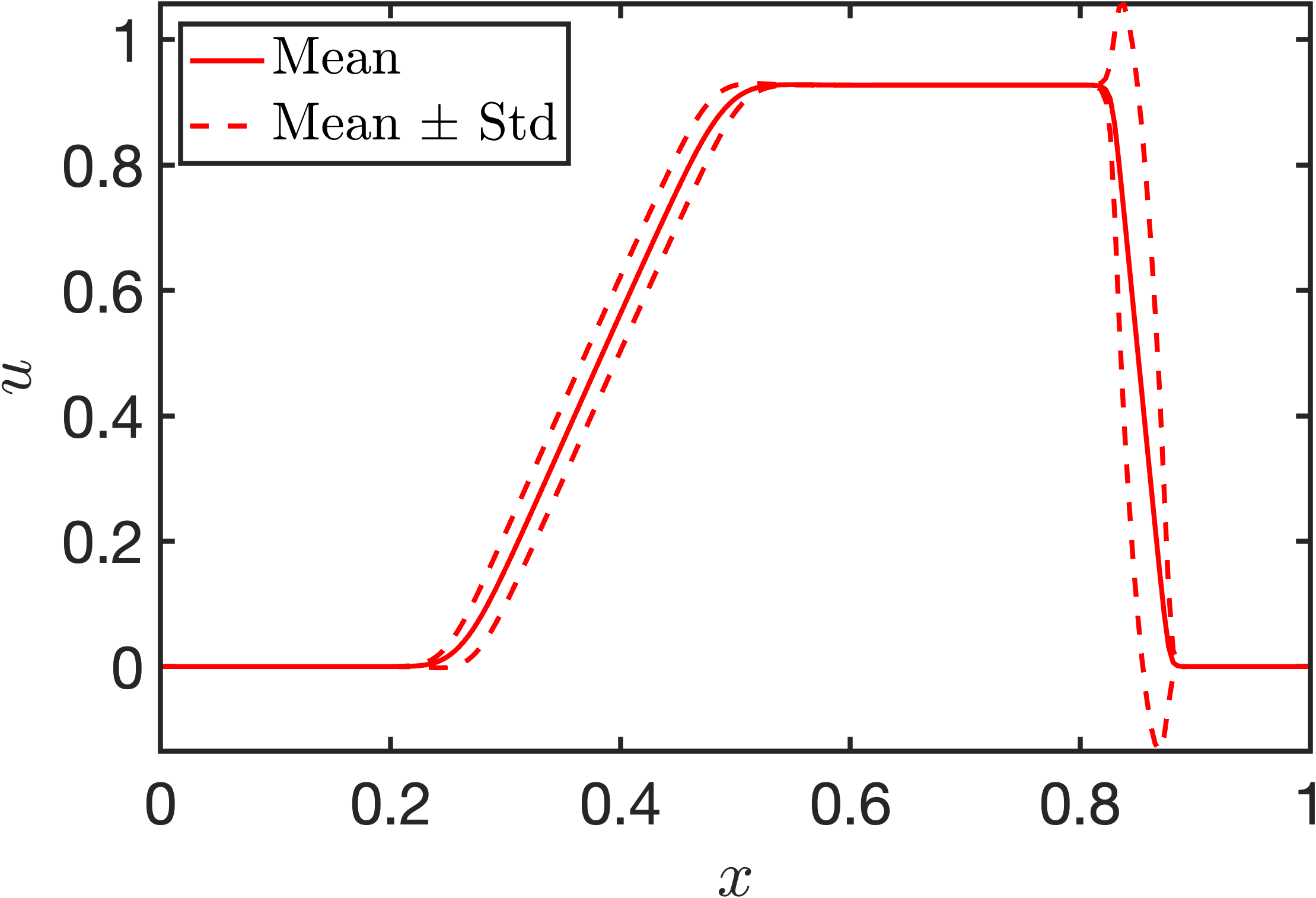}
        \caption{Mean$\pm$std of $u$, $N=13$ ROM}
    \end{subfigure}
    \caption{Stochastic Sod Shock tube 1: solution plots and mean for $N_H = N = 13$ (resolved) hyper-reduced ROM.}
    \label{fig:Sod1_N13}
\end{figure}

\subsection{Sod shock tube with one stochastic variable}
\label{sec:Sod1}

We first consider the Sod shock tube problem \cite{Sod78} over domain $D_x=\LRs{0,1}$ with uncertainty in the initial shock position:
\begin{equation}
    \bm{W}_0(x,y_1) = [\rho_0,u_0,p_0]^T = 
\begin{cases}
    [1.0,0.0,1.0]^T & x<y_1, \\
    [0.125,0.0,0.1]^T & x> y_1,
\end{cases}
\label{eq:Sod1_IC}
\end{equation}
where $y_1\sim\mathcal{U}[0.475,0.525]$. We implement outflow boundary condition and set $\gamma=1.4$. The simulation is run until the final time $T=0.2$. This problem setting was first considered for the SFV method in \cite{Abgrall2017}.

We apply the full SFV method to solve the problem using both WENO reconstruction procedures described in Section \ref{sec:WENO} and report their relative differences in \autoref{tb:Sod1}. In this comparison, we fix $N_x=256$ and vary $N_y$. The solution obtained from the procedure of reconstructed states is taken as a reference (as validated in \cite{Tokareva14,Abgrall2017}) to assess the accuracy of the procedure of reconstructed fluxes. By refining the stochastic grid, we can verify that WENO with reconstructed fluxes achieves similar accuracy and even produce identical solutions (with relative difference around machine precision) after $N_y = 32$. We also compare their runtime per time step. The full SFV method using WENO with reconstructed fluxes consistently achieves faster runtime across all cases.

Now we fix $N_y = 32$. We first examine the singular values associated with the flux snapshots in \autoref{fig:Sod1_svd}, where a sharp decay is observed after $N = 13$. Motivated by this behavior, we vary the number of modes $N$ from 3 to 13 and set $N_H = N$. We then plot the ROM error and the runtime per time step in \autoref{fig:Sod1_er_time}. In general, as $N$ increases, the ROM error decreases. The runtime per time step initially increases, exceeding 5 milliseconds per time step for $7 \leq N \leq 12$. However, at $N = 13$, both the error and runtime drop sharply (ROM error $1.81\times 10^{-10}$) as the reduced basis becoming sufficient to accurately resolve the shocks in the reconstructed flux space, which is consistent with the singular value behavior observed in \autoref{fig:Sod1_svd}.

We now examine the reconstructed solutions. In \autoref{fig:Sod1_N5} and \autoref{fig:Sod1_N13}, we plot the density $\rho$ and velocity $u$ obtained from the ROM, along with their corresponding mean and standard deviation, for $N=5$ and $N=13$, respectively. For $N=5$, noticeable oscillations are present in both the physical and stochastic directions, and negative $\rho$ and $u$ at some cells indicate the ROM is under-resolved and not physically meaningful at those cells. In contrast, for $N=13$, the results exhibit no spurious oscillations, demonstrating that this ROM is stable, accurate, and computationally efficient. Note that we omit the solution plots from the full SFV method, as the ROM with $N=13$ yields a very small error and shows no visual difference.

\begin{table}
\caption{Stochastic Shu--Osher shock tube: relative difference between solutions from full SFV method using WENO with reconstructed states and WENO with reconstructed fluxes, together with their runtime per time step in milliseconds.}
\label{tb:ShuOsher}
\begin{tabular}{lccc}
\hline\noalign{\smallskip}
$N_y$ & Difference & Runtime per step (state) & Runtime per step (flux) \\
\noalign{\smallskip}\hline\noalign{\smallskip}
$4$  & $5.23\times 10^{-5}$ & $1.55$  & $1.45$ \\
$8$  & $4.97\times 10^{-5}$ & $2.22$  & $2.07$ \\
$16$ & $1.40\times 10^{-5}$ & $3.93$  & $3.59$ \\
$32$ & $3.65\times 10^{-15}$ & $1.43 \times 10^1$ & $1.40 \times 10^1$ \\
\noalign{\smallskip}\hline
\end{tabular}
\end{table}

\begin{figure}
    \includegraphics[width=0.5\textwidth]{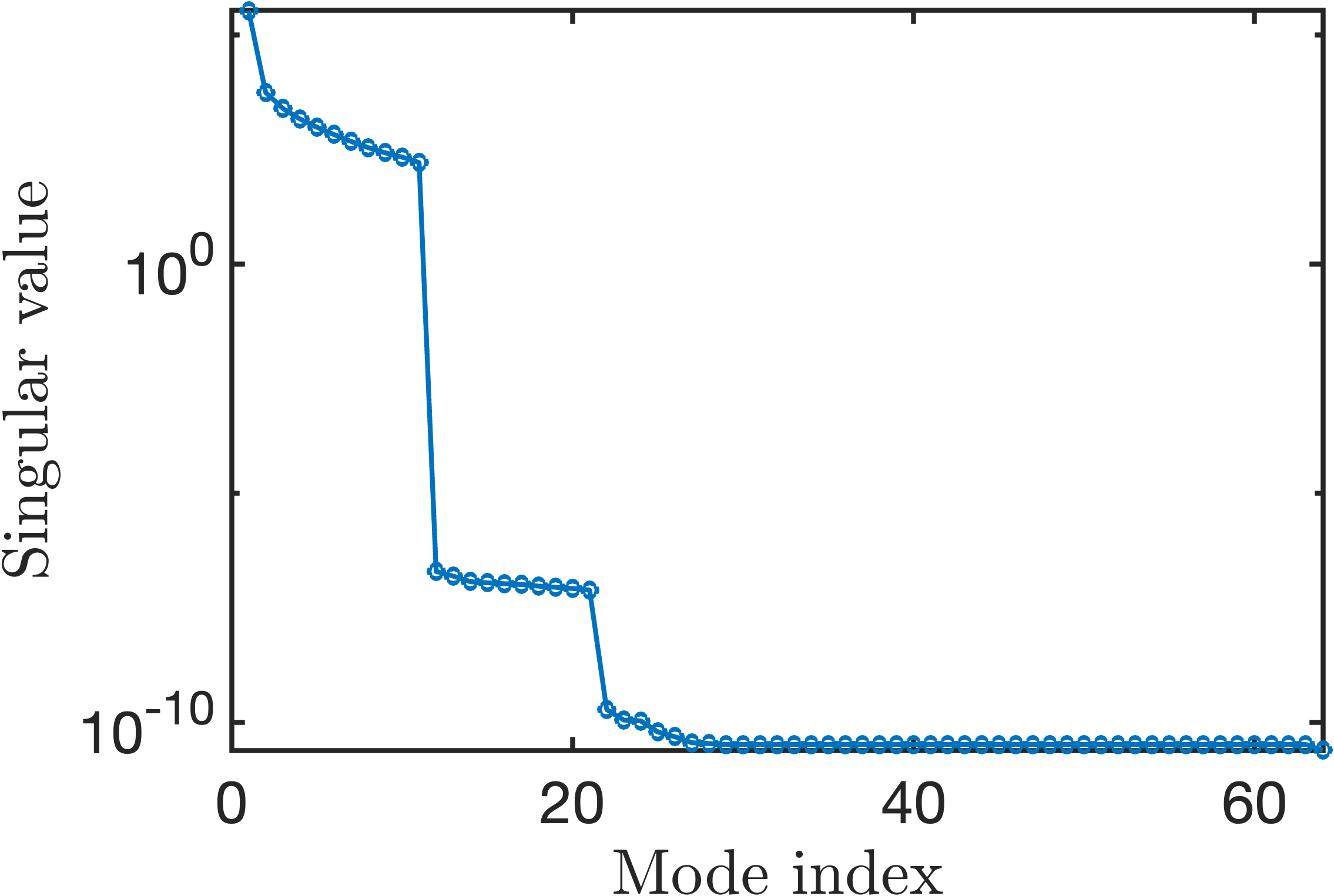}
    \caption{Stochastic Shu--Osher shock tube: singular values of the flux snapshot matrix.}
    \label{fig:ShuOsher_svd}
\end{figure}

\begin{figure}
    \centering
    \begin{subfigure}{0.49\textwidth}
        \centering
        \includegraphics[width=\textwidth]{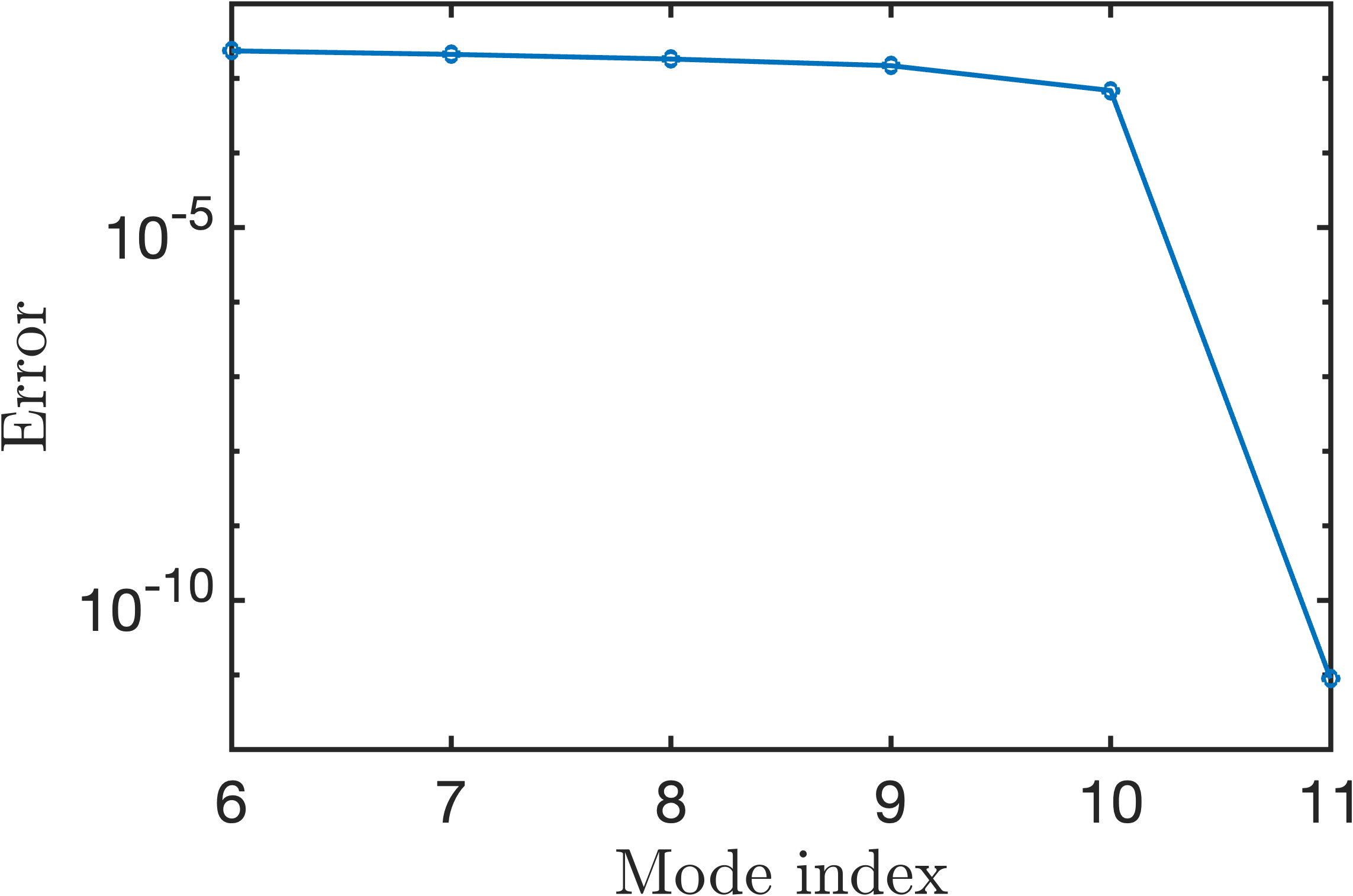}
        \caption{ROM error plot}
    \end{subfigure}
    \begin{subfigure}{0.49\textwidth}
        \centering
        \includegraphics[width=\textwidth]{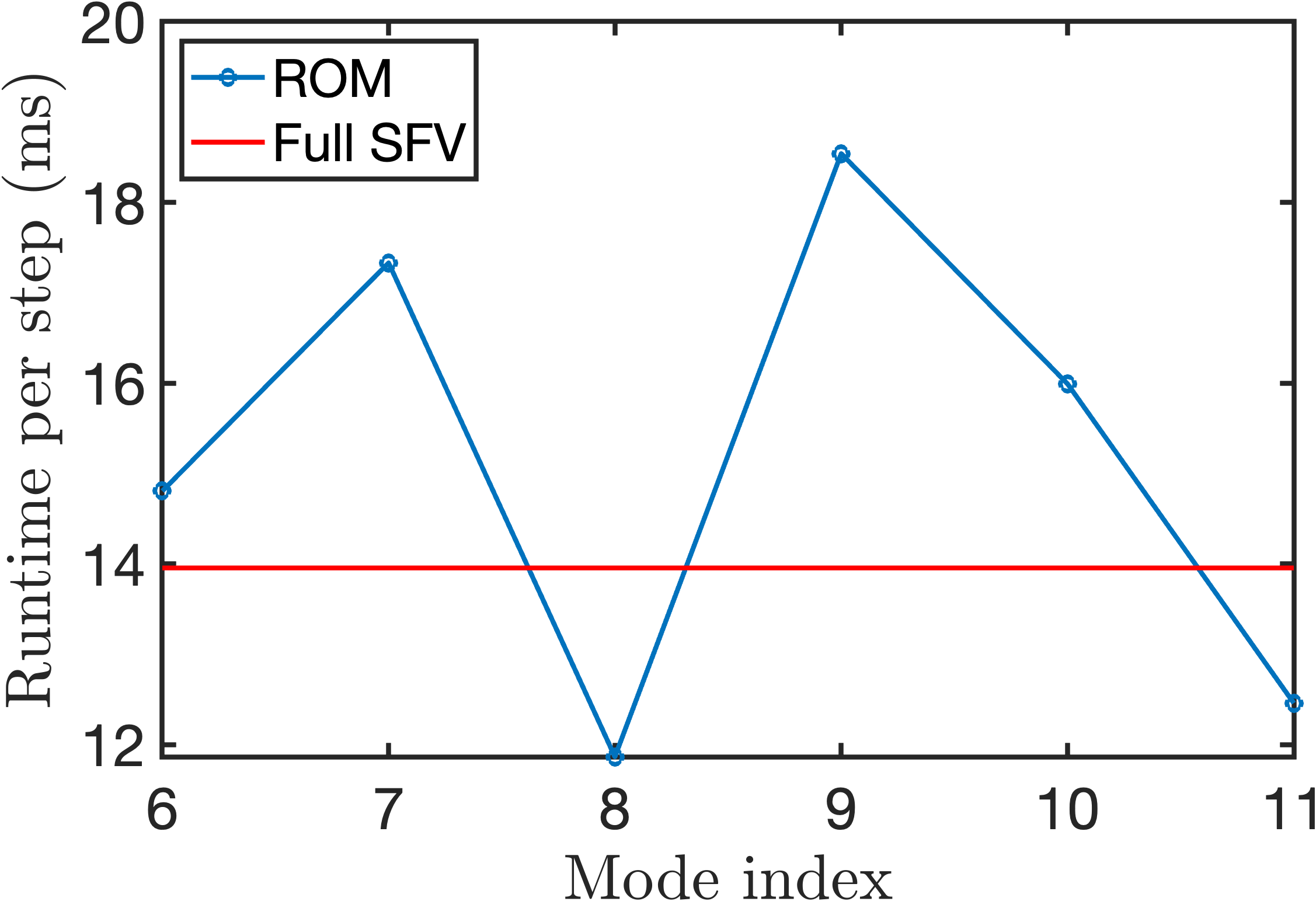}
        \caption{ROM time-stepping runtime plot}
    \end{subfigure}
    \caption{Stochcastic Shu--Osher shock tube: error and runtime per time step for $N_H = N$ hyper-reduced ROMs.}
    \label{fig:ShuOsher_er_time}
\end{figure}

\begin{figure}
    \centering
    \begin{subfigure}{0.49\textwidth}
        \centering
        \includegraphics[width=\textwidth]{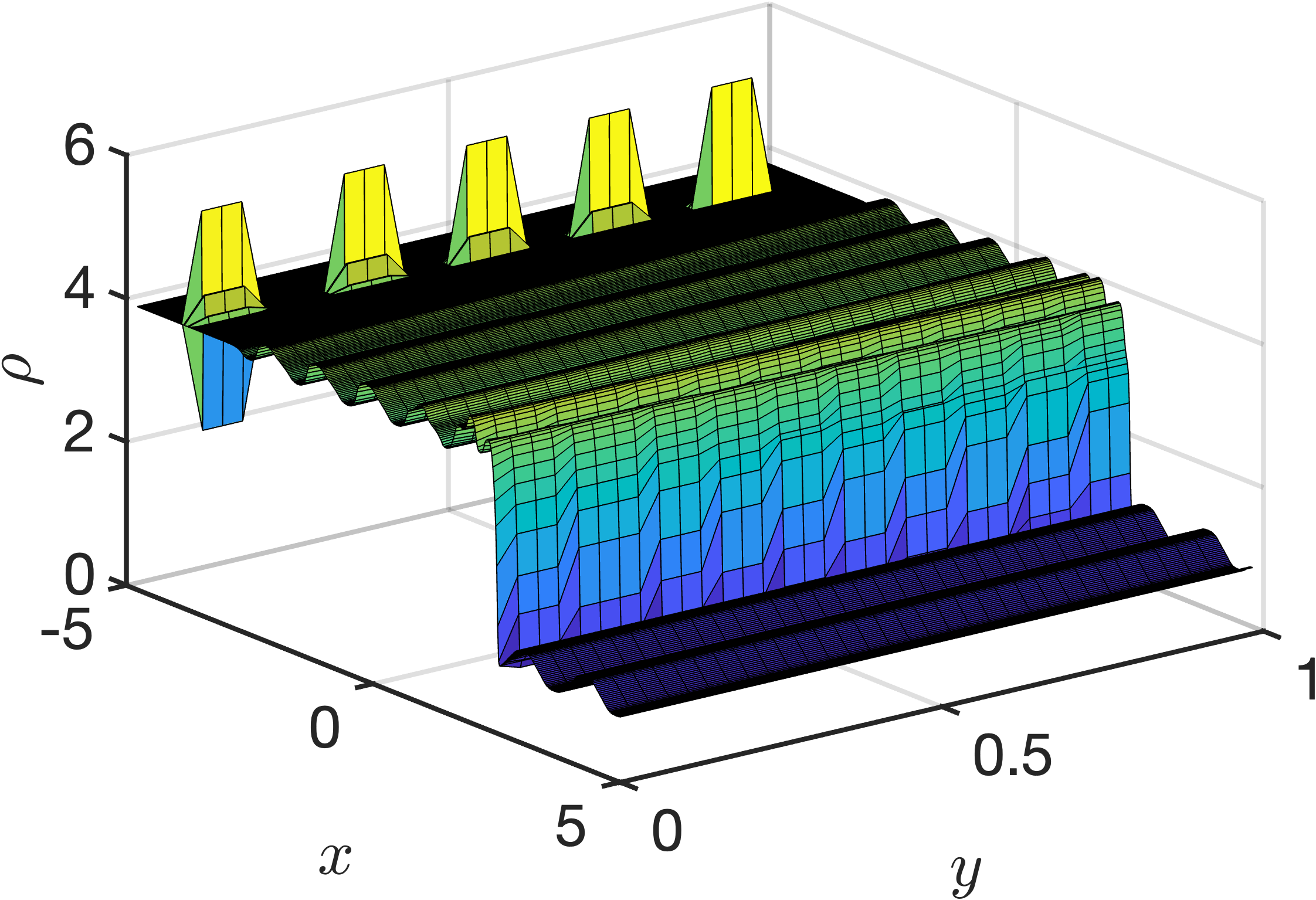}
        \caption{$\rho$ plot, $N=6$ ROM}
    \end{subfigure}
    \begin{subfigure}{0.49\textwidth}
        \centering
        \includegraphics[width=\textwidth]{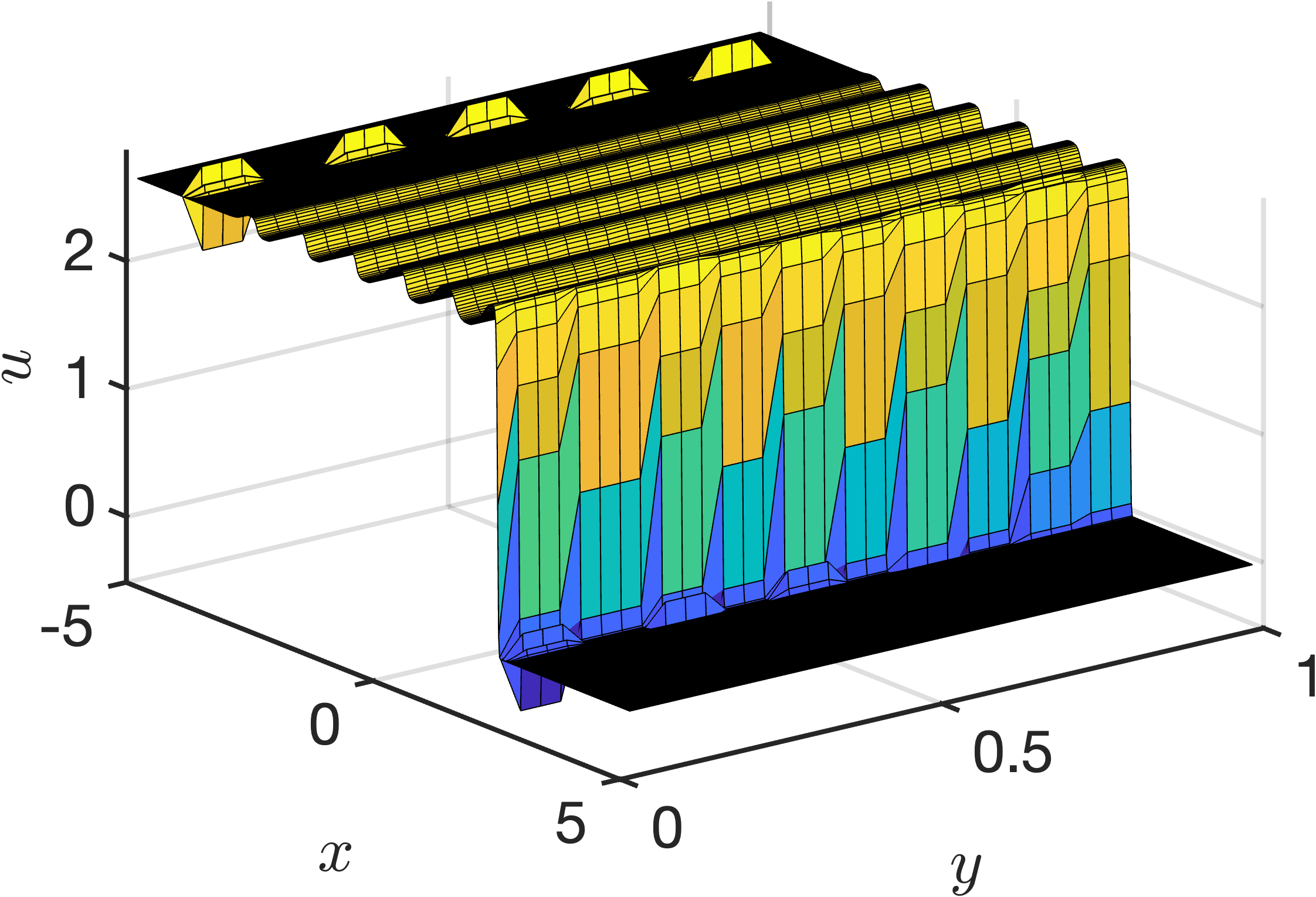}
        \caption{$u$ plot, $N=6$ ROM}
    \end{subfigure}
     \begin{subfigure}{0.49\textwidth}
        \centering
        \includegraphics[width=\textwidth]{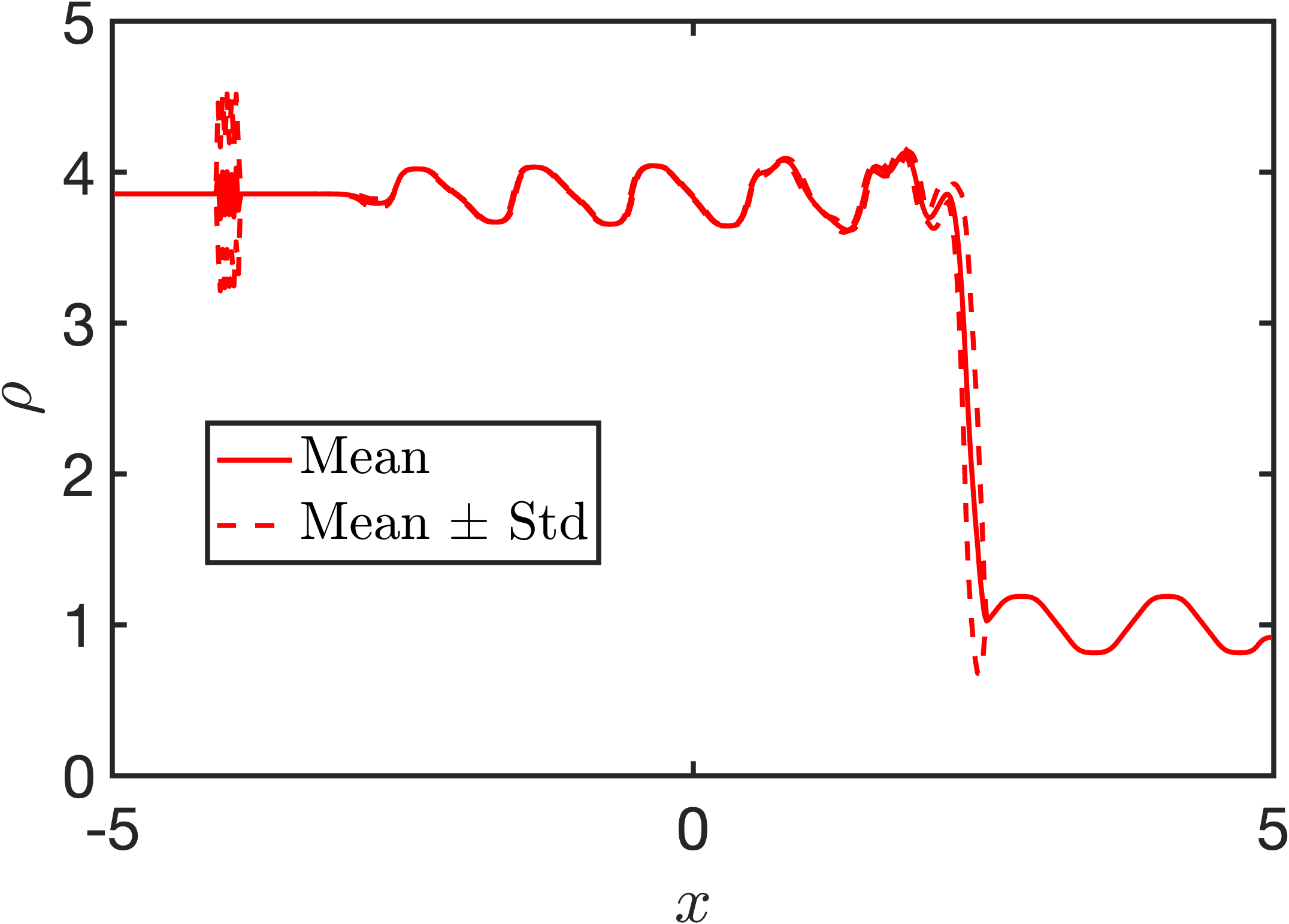}
        \caption{Mean$\pm$std of $\rho$, $N=6$ ROM}
    \end{subfigure}
    \begin{subfigure}{0.49\textwidth}
        \centering
        \includegraphics[width=\textwidth]{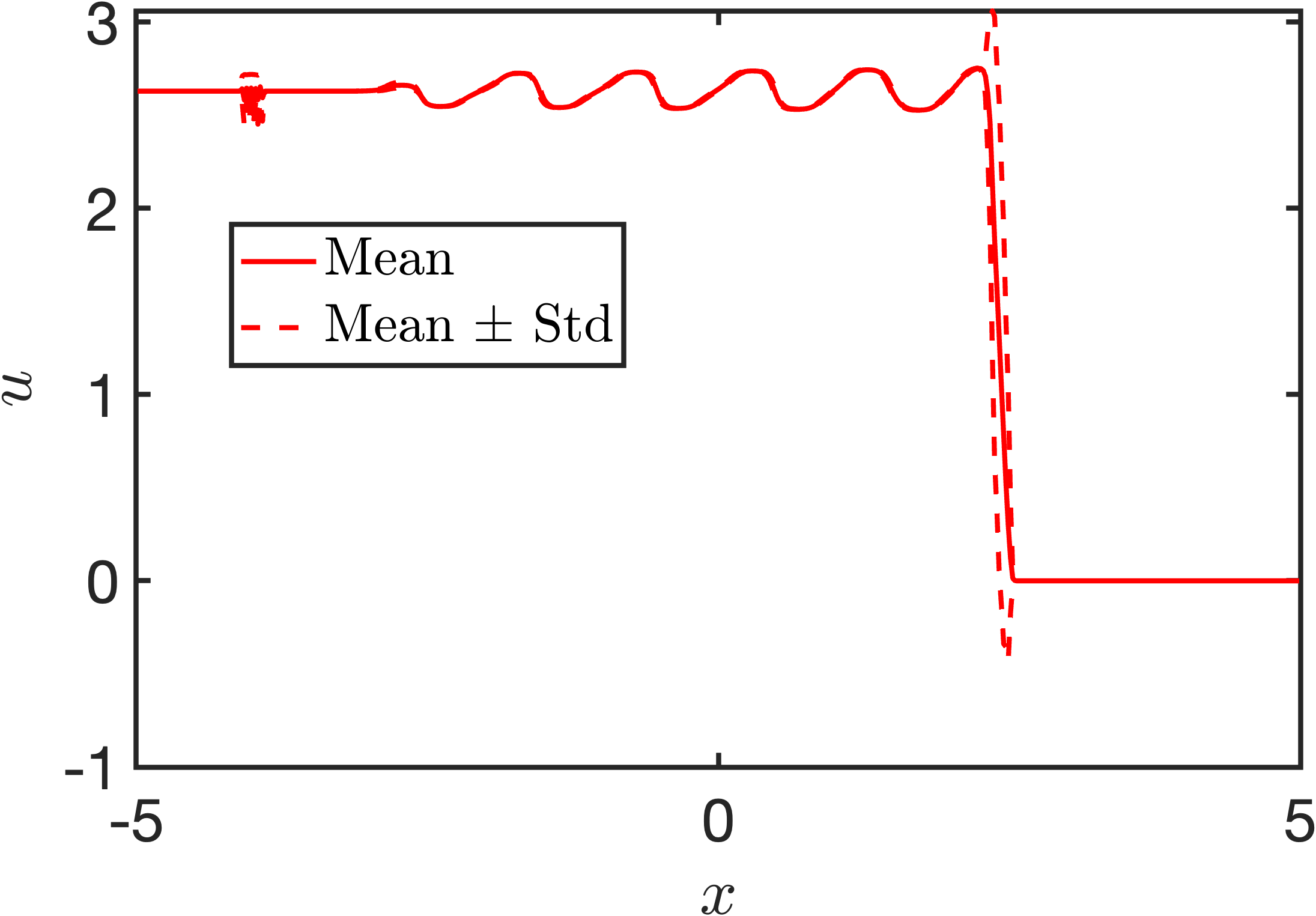}
        \caption{Mean$\pm$std of $u$, $N=6$ ROM}
    \end{subfigure}
    \caption{Stochastic Shu--Osher shock tube: solution plots and mean for $N_H = N = 6$ (under-resolved) hyper-reduced ROM.}
    \label{fig:ShuOsher_N6}
\end{figure}

\begin{figure}
    \centering
    \begin{subfigure}{0.49\textwidth}
        \centering
        \includegraphics[width=\textwidth]{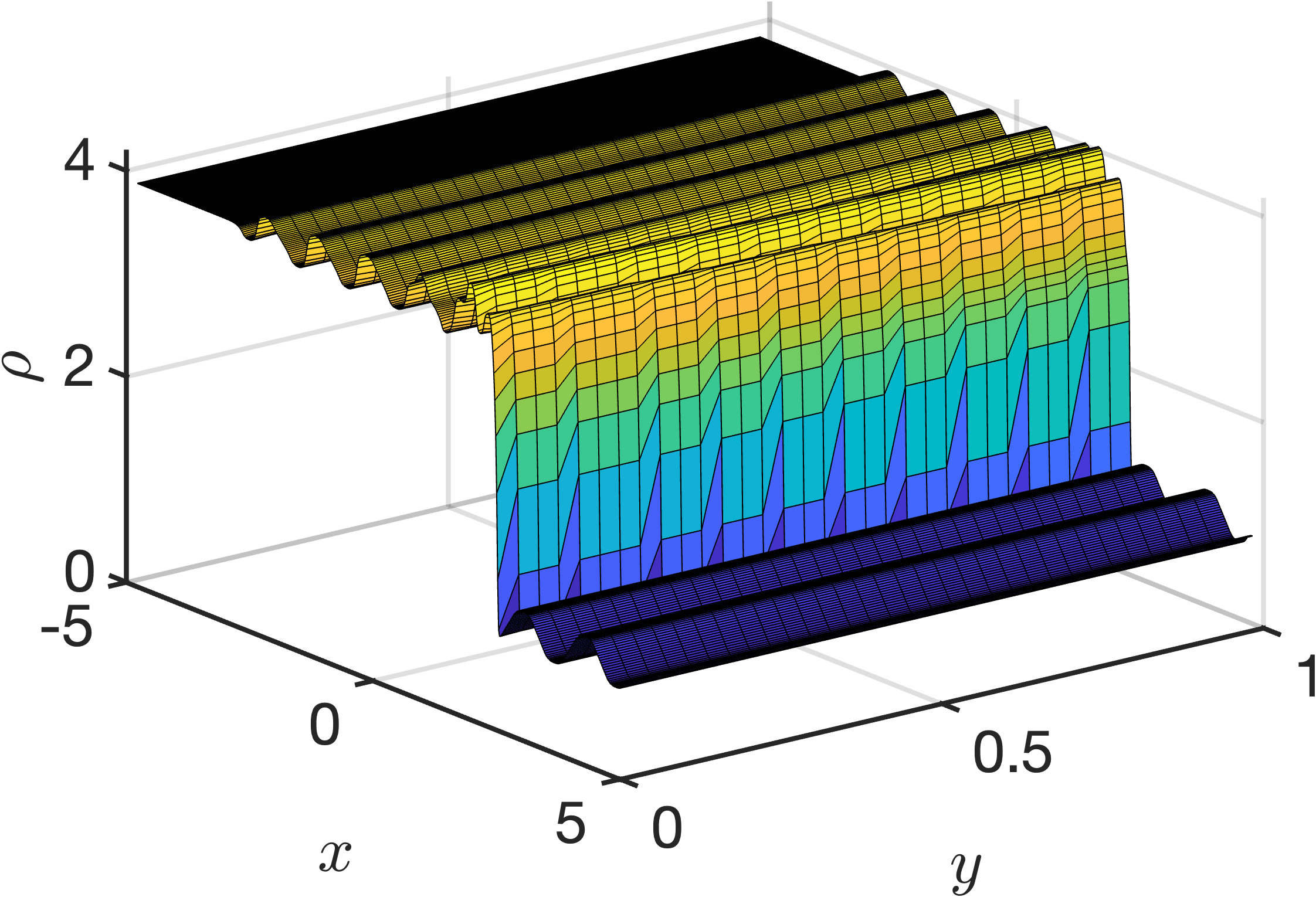}
        \caption{$\rho$ plot, $N=11$ ROM}
    \end{subfigure}
    \begin{subfigure}{0.49\textwidth}
        \centering
        \includegraphics[width=\textwidth]{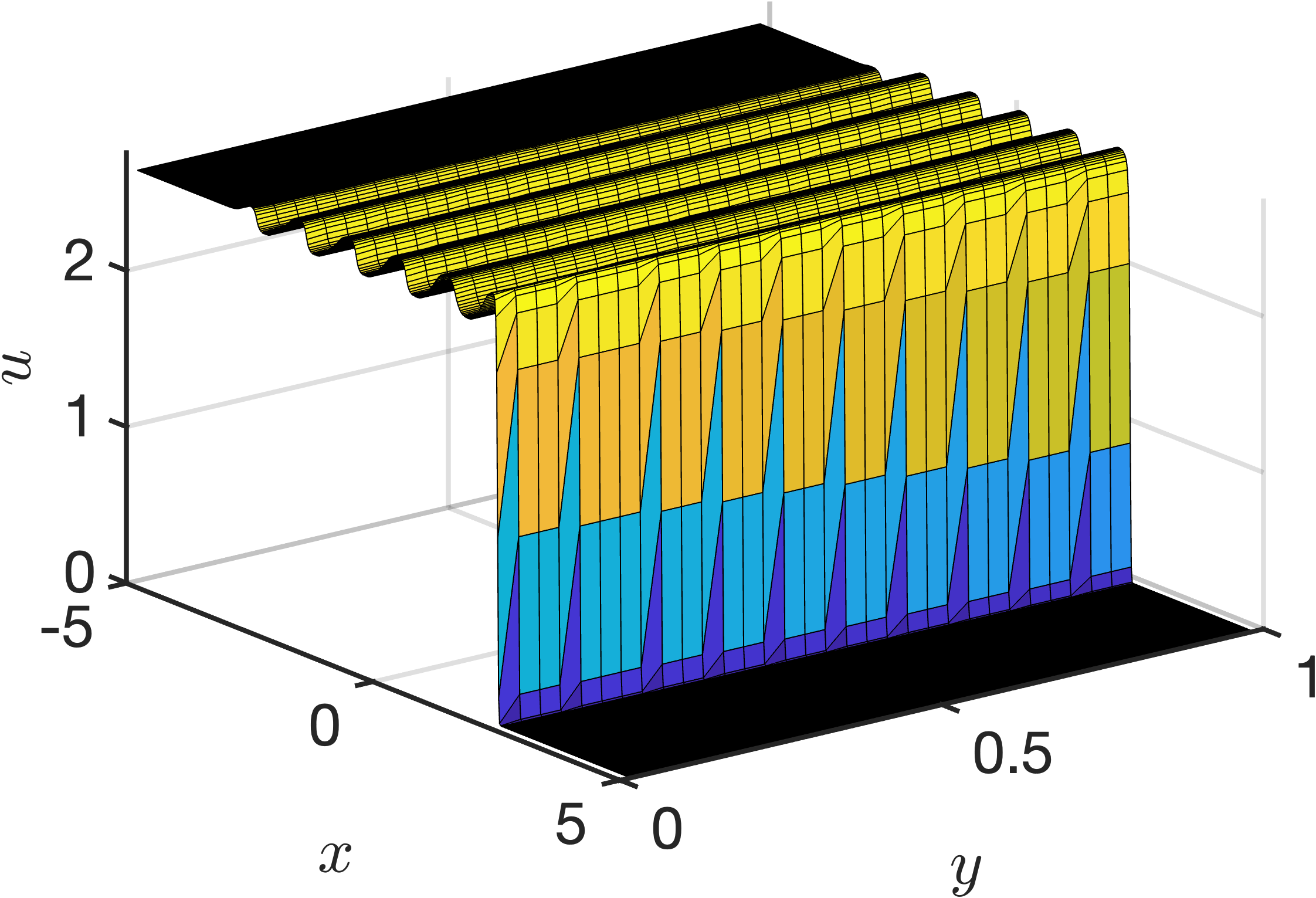}
        \caption{$u$ plot, $N=11$ ROM}
    \end{subfigure}
     \begin{subfigure}{0.49\textwidth}
        \centering
        \includegraphics[width=\textwidth]{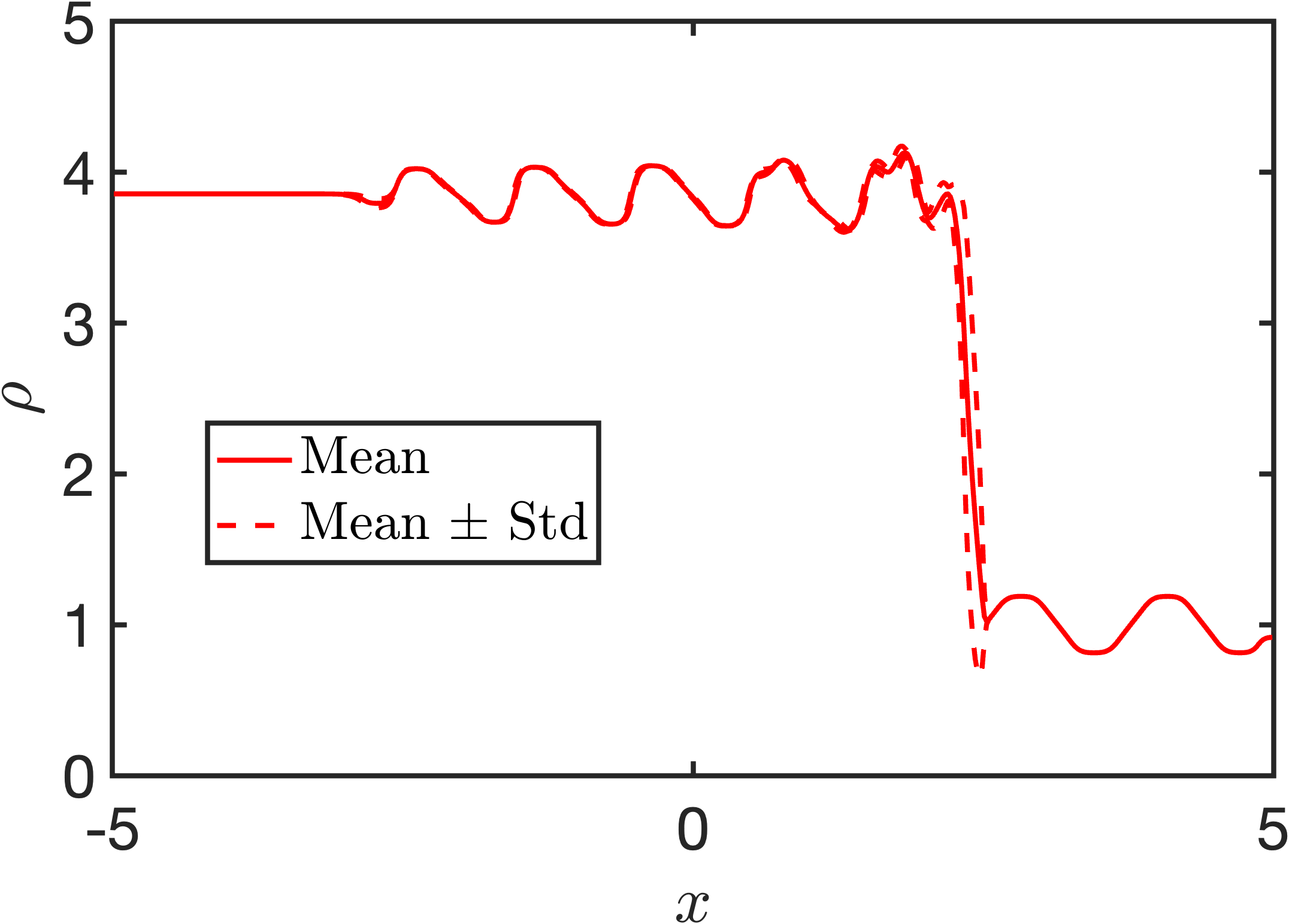}
        \caption{Mean$\pm$std of $\rho$, $N=11$ ROM}
    \end{subfigure}
    \begin{subfigure}{0.49\textwidth}
        \centering
        \includegraphics[width=\textwidth]{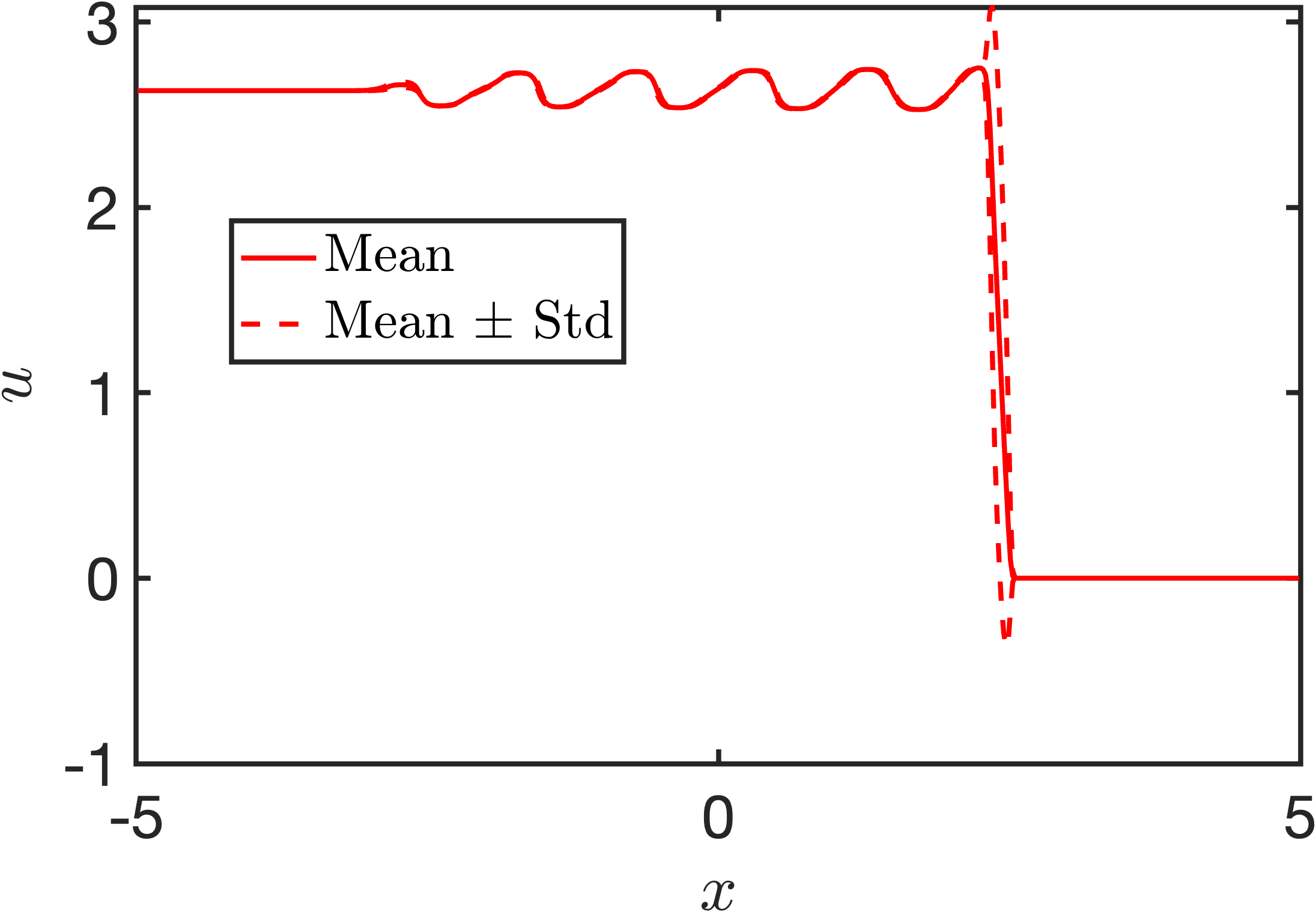}
        \caption{Mean$\pm$std of $u$, $N=11$ ROM}
    \end{subfigure}
    \caption{Stochastic Shu--Osher shock tube: solution plots and mean for $N_H = N = 11$ (resolved) hyper-reduced ROM.}
    \label{fig:ShuOsher_N11}
\end{figure}

\subsection{Shu--Osher shock tube  with one stochastic variable}
We next consider the Shu--Osher shock tube problem \cite{SHU1989} over the domain $D_x = [-5,5]$ with uncertainty in the shock location:
\[
\bm{W}_0(x,y) = [\rho_0,u_0,p_0]^T =
\begin{cases}
    [3.857143,\,2.629369,\,10.3333]^T & x < y_1, \\
    \big[1 + 0.2\sin(5x),\,0,\,1\big]^T & x > y_1,
\end{cases}
\]
where $y_1 \sim \mathcal{U}[-4.1,-3.9]$. We implement outflow boundary condition and set $\gamma=1.4$. The simulation is run until the final time $T = 1.8$.

We again first compare the two WENO reconstruction schemes for the full SFV method. Fixing $N_x = 512$  and varying $N_y$, we report their relative difference and timing in \autoref{tb:ShuOsher}. We observe that for this problem, the relative difference stay around the same small level and reach close to machine precision at $N_y = 32$, and WENO with reconstructed fluxes is consistently faster. 

We now fix $N_y = 32$. We first examine the singular values of the flux snapshots in \autoref{fig:ShuOsher_svd}, where a sharp decay is observed after $N = 11$. Motivated by this behavior, we vary the number of modes $N$ from 6 to 11 and report the hyper-reduced ROM error and runtime per time step with $N_H = N$ in \autoref{fig:ShuOsher_er_time}. We observe behavior similar to the previous example: the ROM error decreases consistently as $N$ increases, while the runtime may even exceed that of the full SFV method with reconstructed fluxes. At $N = 11$, both the error and runtime drop sharply (ROM error $8.87\times 10^{-12}$) since the first 11 modes are sufficient to resolve the shocks in the reconstructed flux space, as indicated by the decay of the singular values.

Finally, we plot the density $\rho$ and velocity $u$ obtained from the ROM, along with their corresponding mean and standard deviation, for $N=6$ in \autoref{fig:ShuOsher_N6} and $N=11$ in \autoref{fig:ShuOsher_N11}. At $N=6$, we observe oscillations similar to those in the previous examples in both $\rho$ and $u$, which are also reflected in the corresponding standard deviations. At $N=11$, the solutions appear smooth and free of spurious oscillations. As before, we do not include the solution plots from the full SFV method, since the $N=11$ ROM error is small.

\begin{table}
\caption{Stochastic Sod shock tube 2: relative difference between solutions from full SFV method using WENO with reconstructed states and WENO with reconstructed fluxes, together with their runtime per time step in milliseconds.}
\label{tb:Sod2}
\begin{tabular}{lccc}
\hline\noalign{\smallskip}
$N_y$ & Difference & Runtime per step (state) & Runtime per step (flux) \\
\noalign{\smallskip}\hline\noalign{\smallskip}
$4^2$  & $9.33\times 10^{-4}$ & $1.04 \times 10^1$  & $8.70$ \\
$8^2$  & $7.35\times 10^{-4}$ & $1.65 \times 10^1$  & $1.53\times 10^1$ \\
$16^2$ & $9.21\times 10^{-6}$ & $4.18 \times 10^1$  & $3.63 \times 10^1$ \\
$32^2$ & $2.61\times 10^{-6}$ & $2.06 \times 10^2$ & $1.83 \times 10^2$ \\
\noalign{\smallskip}\hline
\end{tabular}
\end{table}

\begin{figure}
       \centering
    \begin{subfigure}{0.49\textwidth}
        \centering
        \includegraphics[width=\textwidth]{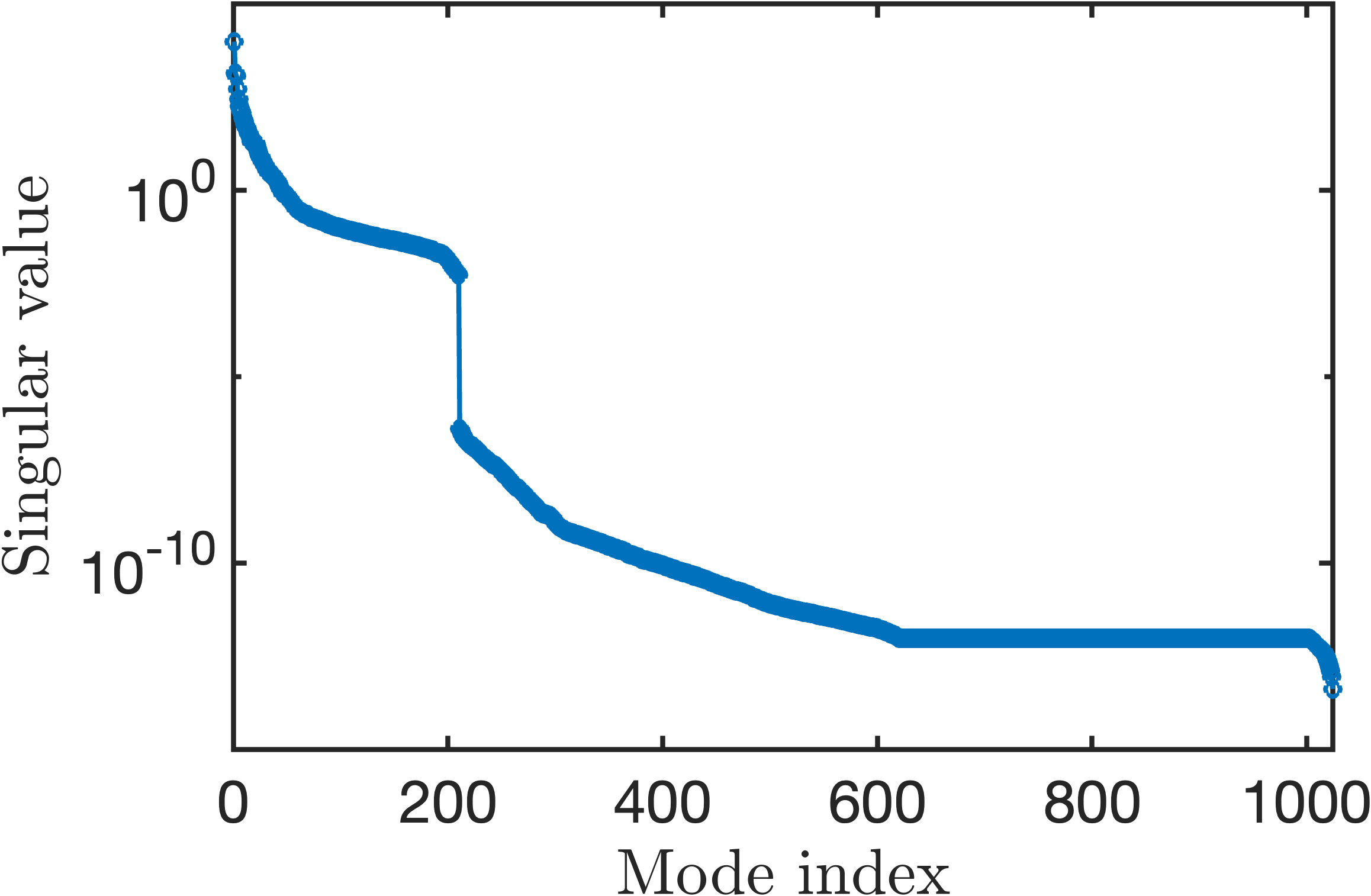}
    \end{subfigure}
    \begin{subfigure}{0.49\textwidth}
        \centering
        \includegraphics[width=\textwidth]{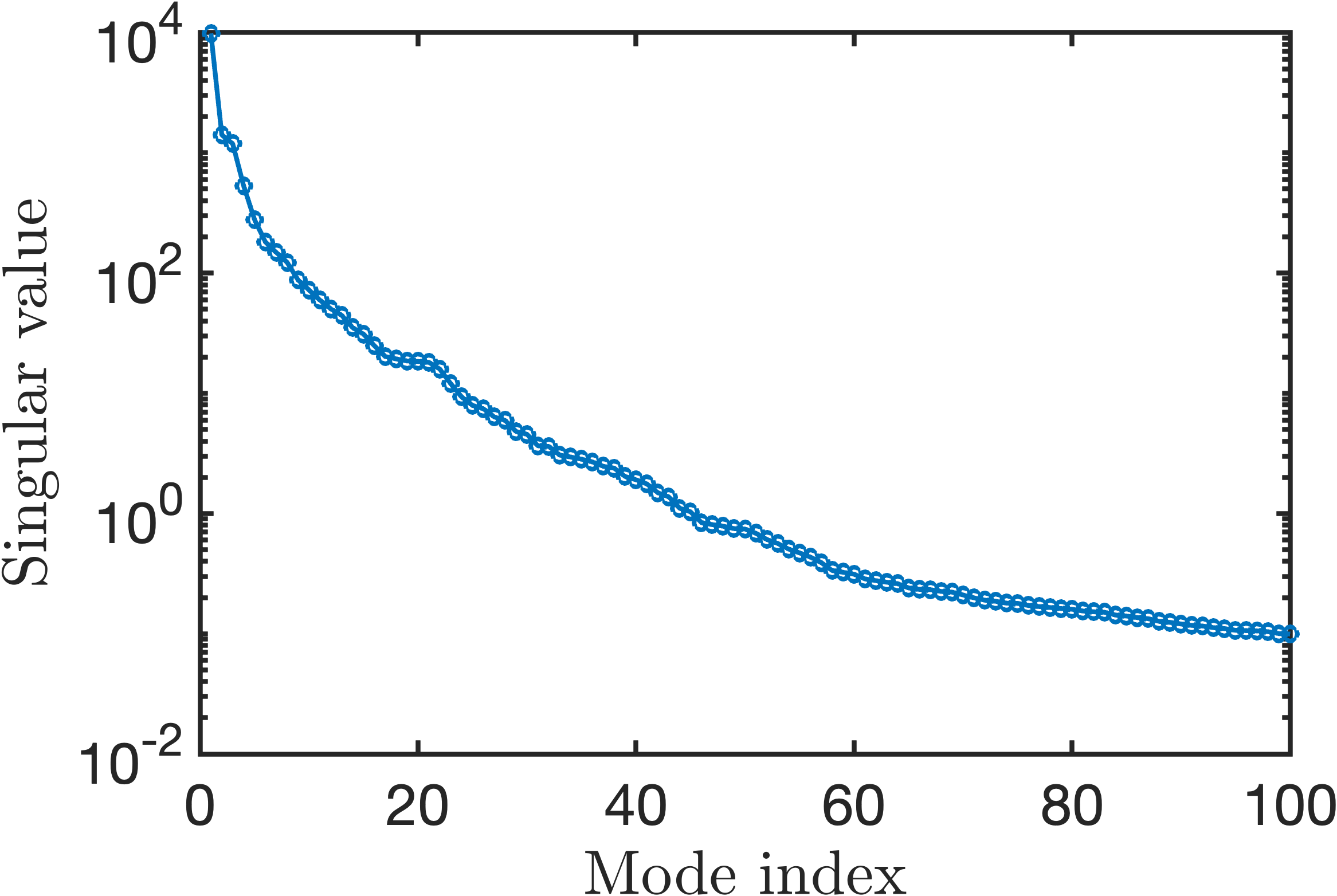}
    \end{subfigure}
    \caption{Stochastic Sod shock tube 2: singular values of the flux snapshot matrix with a zoomed view for the first 100 modes.}
    \label{fig:Sod2_svd}
\end{figure}

\begin{figure}
    \centering
    \begin{subfigure}{0.49\textwidth}
        \centering
        \includegraphics[width=\textwidth]{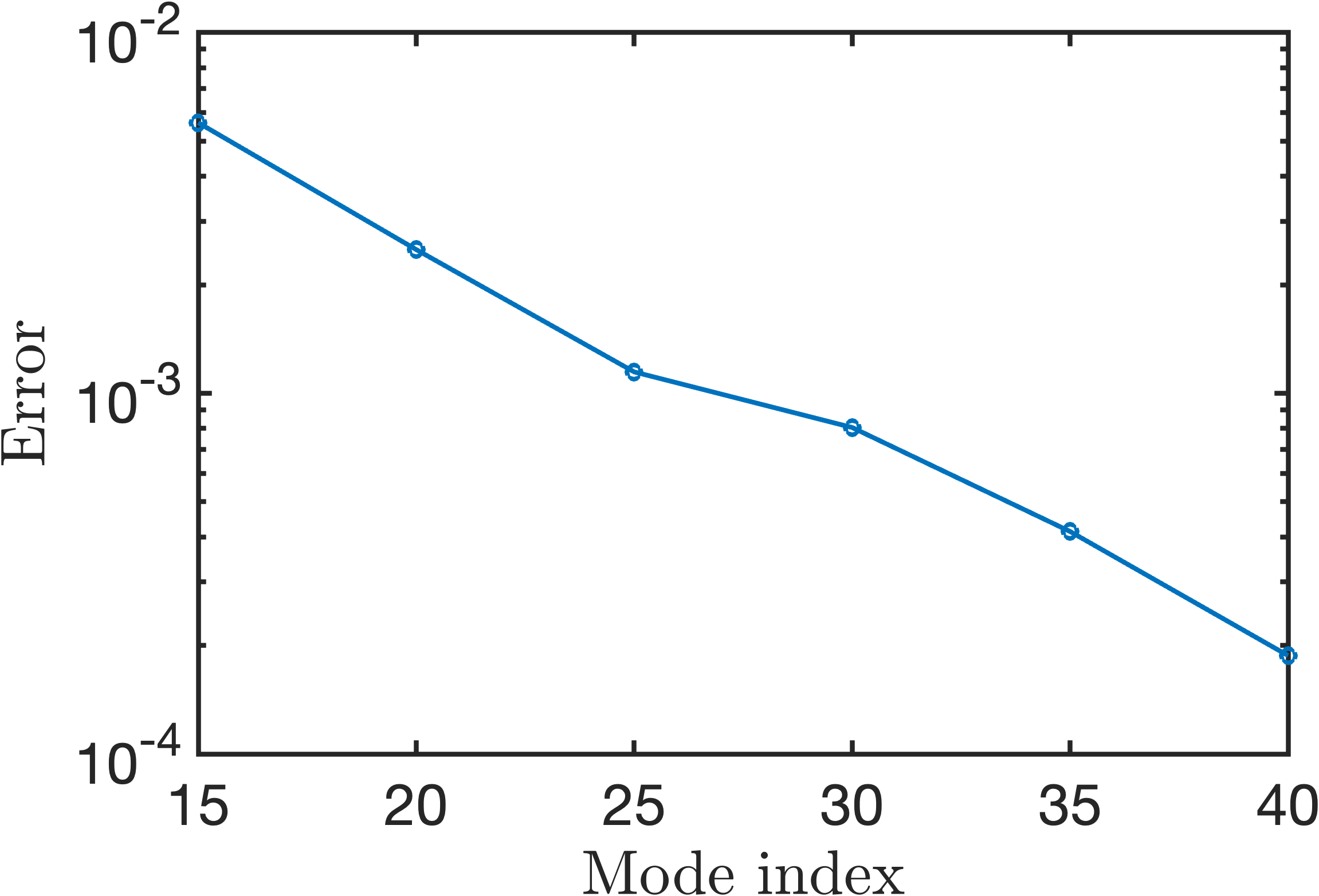}
        \caption{ROM error plot}
    \end{subfigure}
    \begin{subfigure}{0.49\textwidth}
        \centering
        \includegraphics[width=\textwidth]{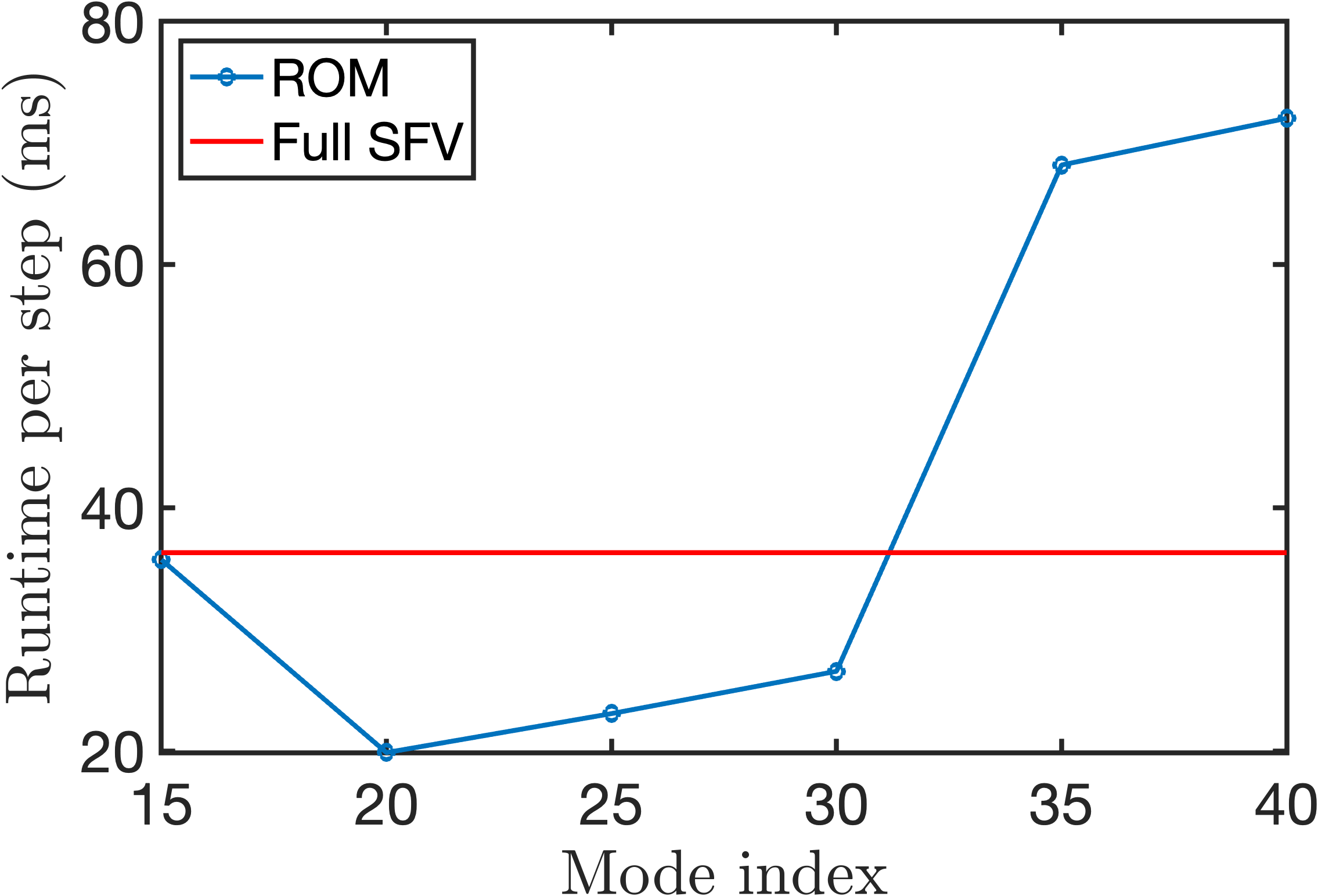}
        \caption{ROM time-stepping runtime plot}
    \end{subfigure}
    \caption{Stochastic Sod shock tube 2: error and runtime per time step for $N_H = N$ hyper-reduced ROMs.}
    \label{fig:Sod2_er_time}
\end{figure}

\begin{figure}
    \centering
    \begin{subfigure}{0.49\textwidth}
        \centering
        \includegraphics[width=\textwidth]{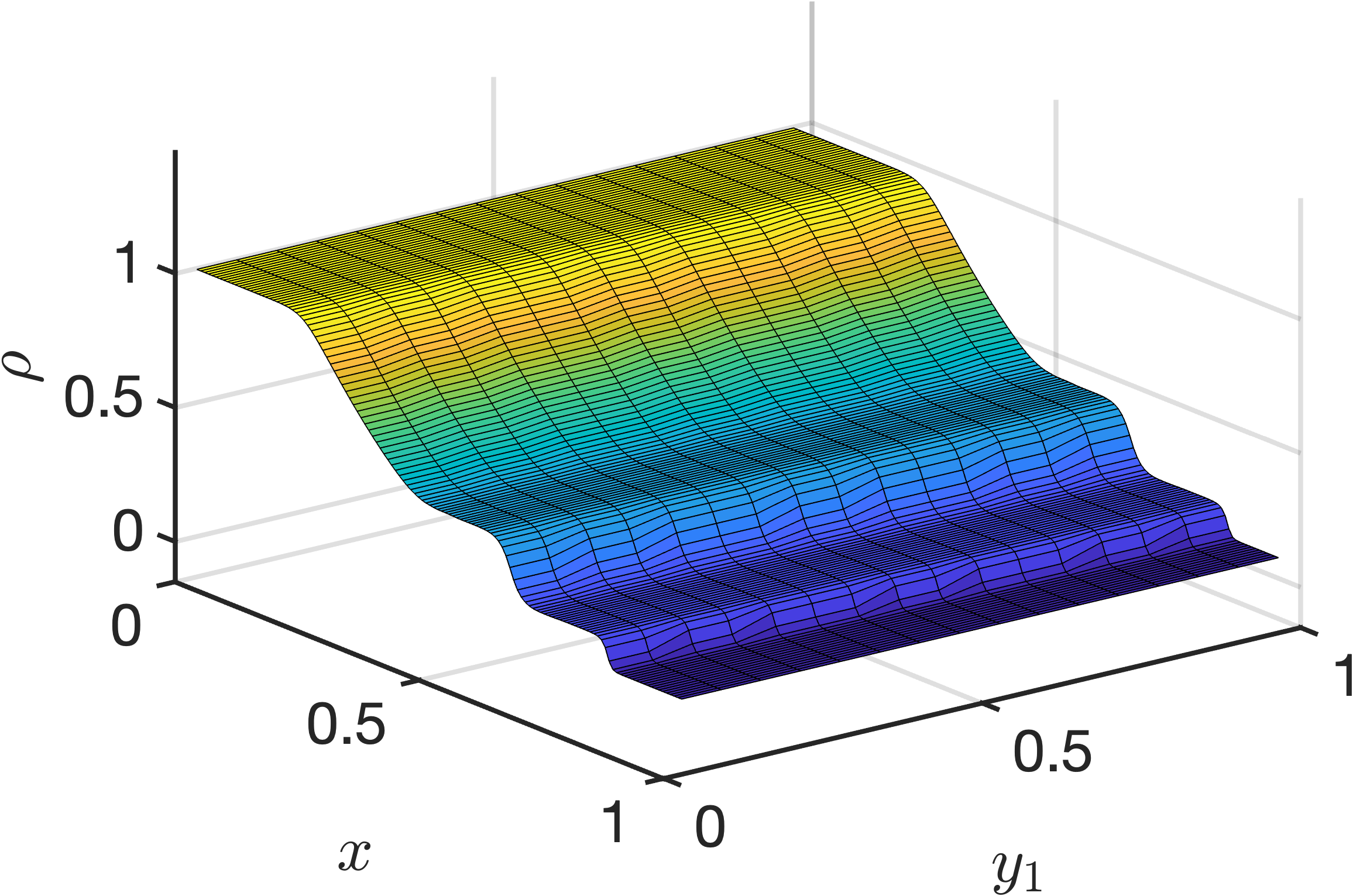}
        \caption{$\rho$ plot at $y_2 =1$}
    \end{subfigure}
    \begin{subfigure}{0.49\textwidth}
        \centering
        \includegraphics[width=\textwidth]{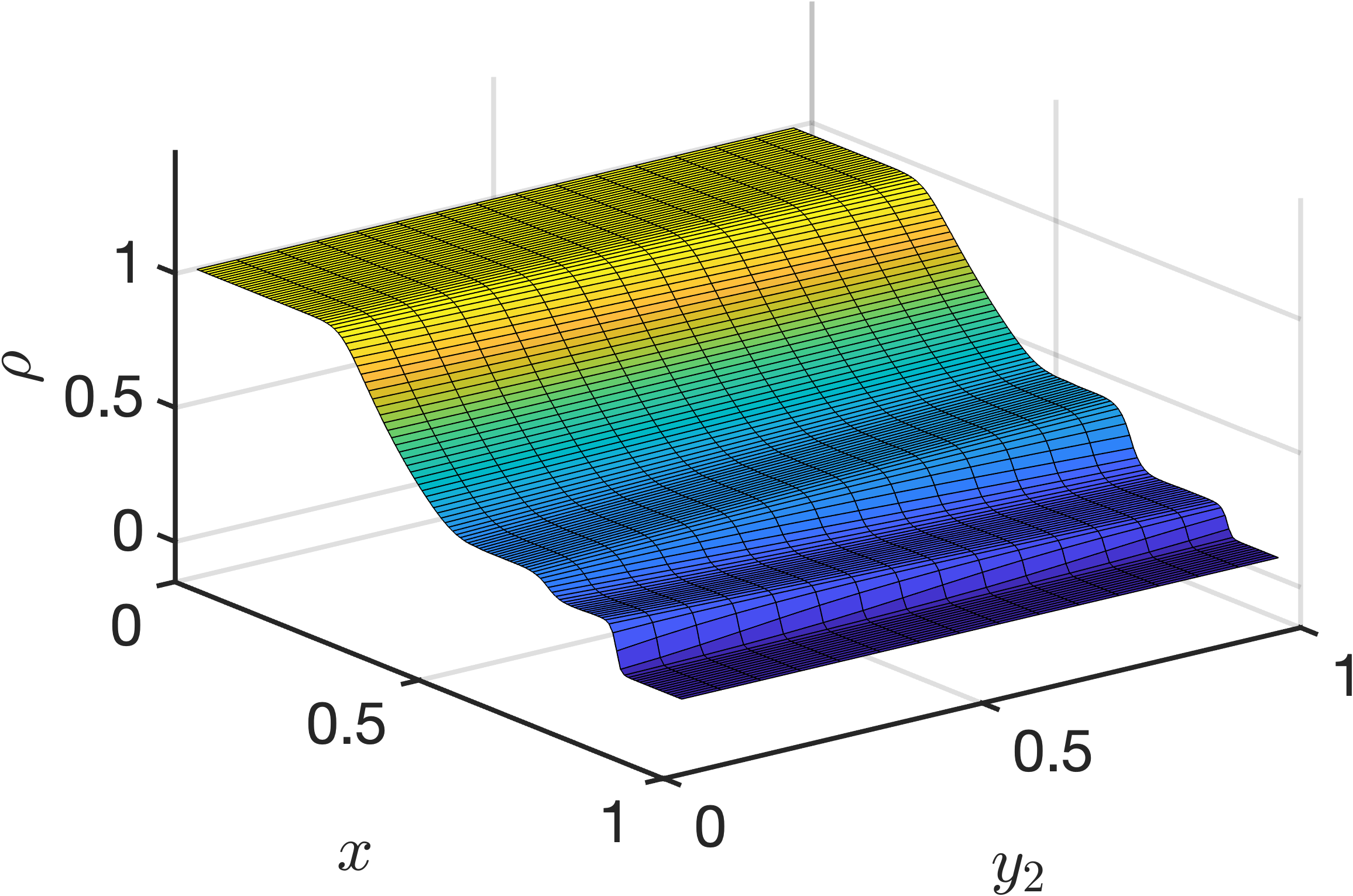}
        \caption{$\rho$ plot at $y_1 =1$}
    \end{subfigure}
     \begin{subfigure}{0.49\textwidth}
        \centering
        \includegraphics[width=\textwidth]{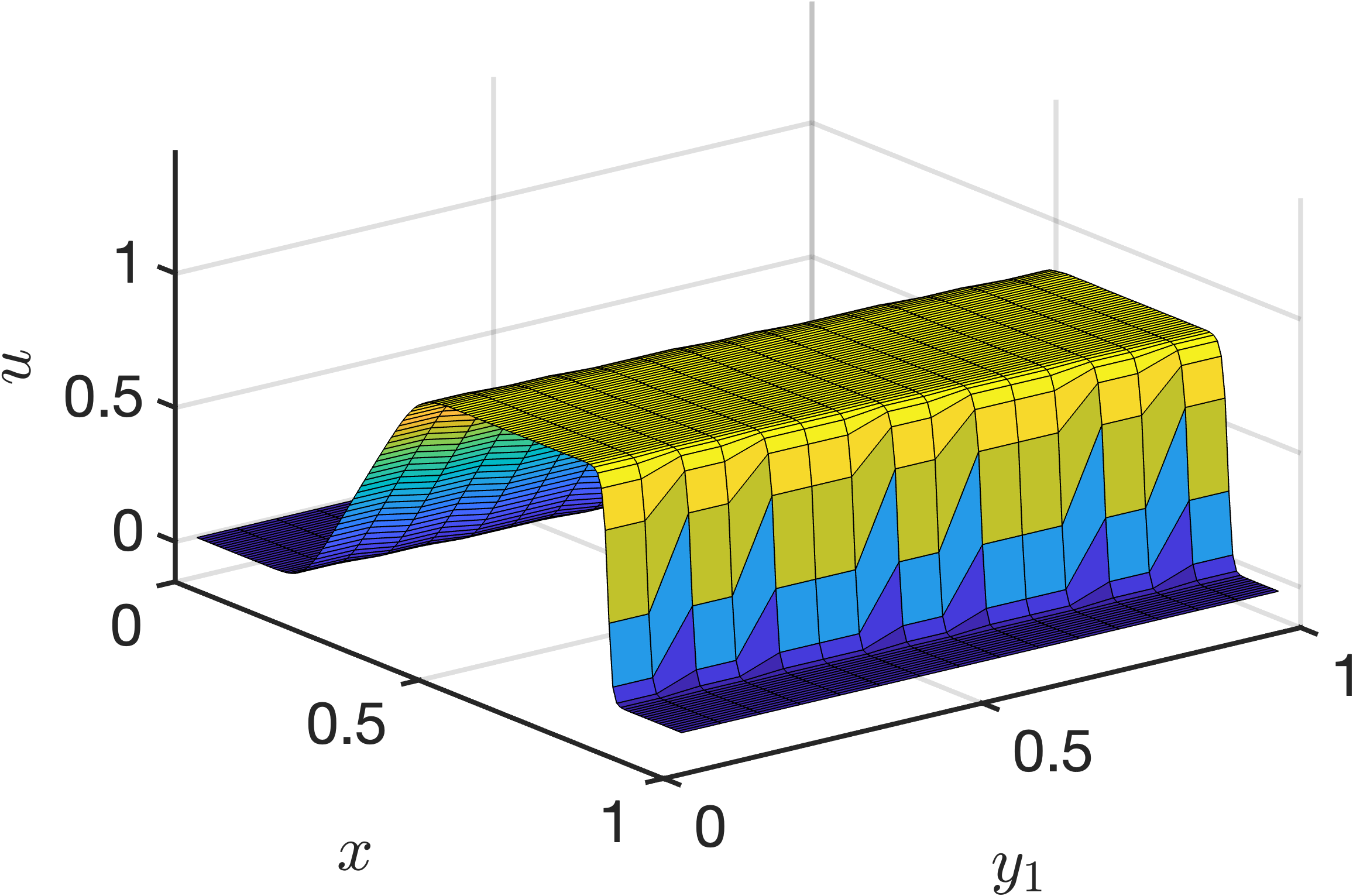}
        \caption{$u$ plot at $y_2 =1$}
    \end{subfigure}
    \begin{subfigure}{0.49\textwidth}
        \centering
        \includegraphics[width=\textwidth]{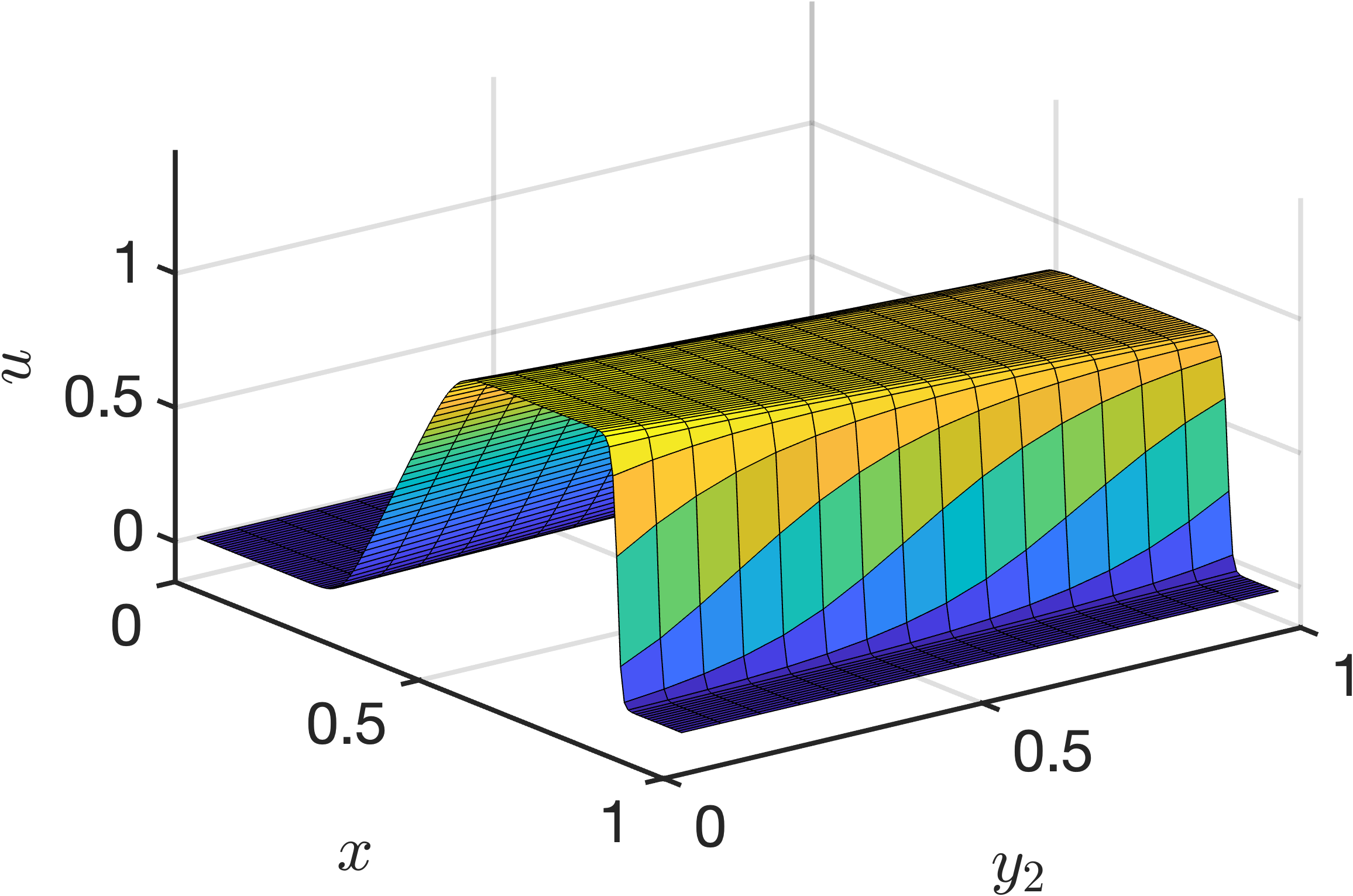}
        \caption{$u$ plot at $y_1 =1$}
    \end{subfigure}
     \begin{subfigure}{0.49\textwidth}
        \centering
        \includegraphics[width=\textwidth]{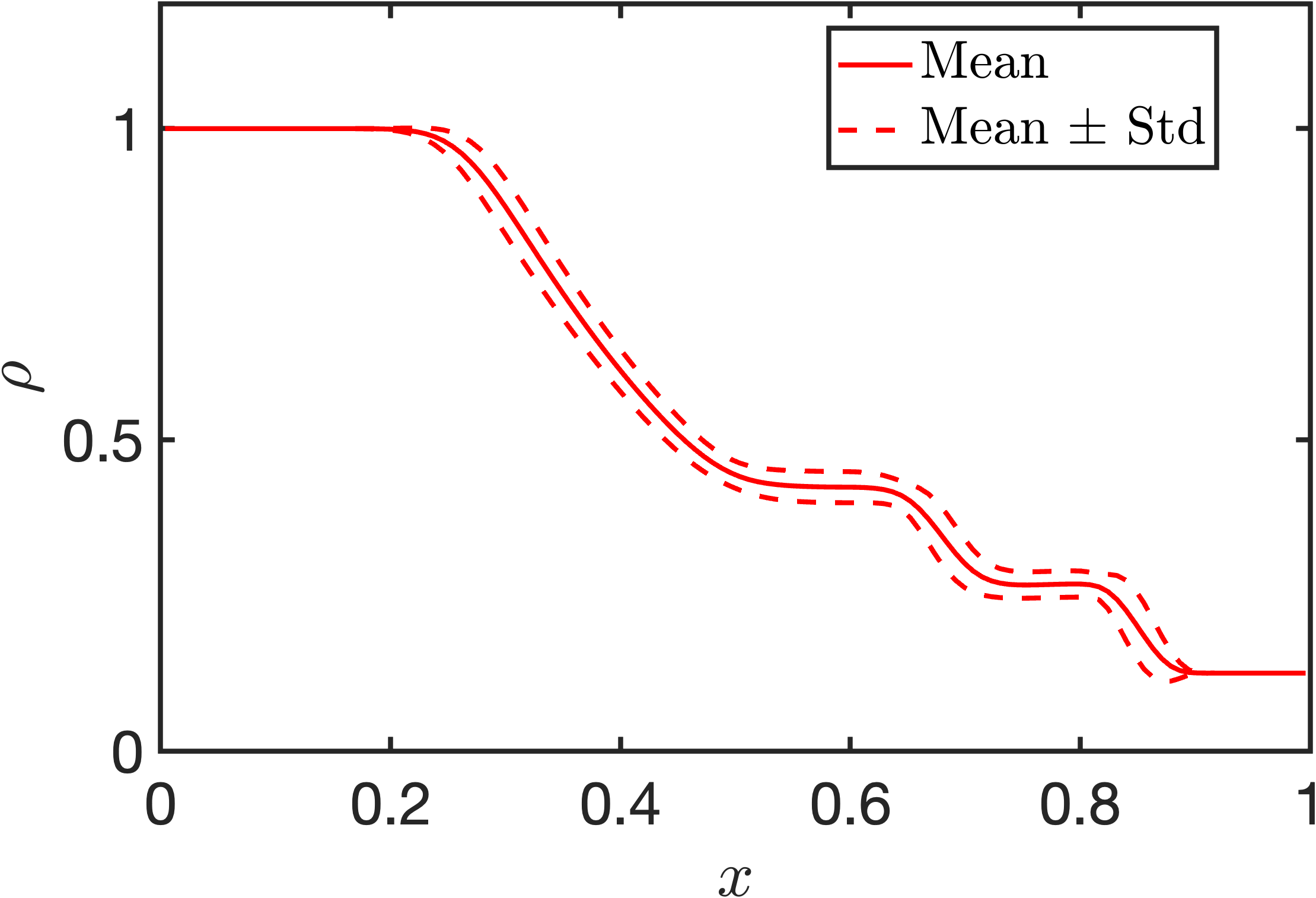}
        \caption{Mean$\pm$std of $\rho$}
    \end{subfigure}
    \begin{subfigure}{0.49\textwidth}
        \centering
        \includegraphics[width=\textwidth]{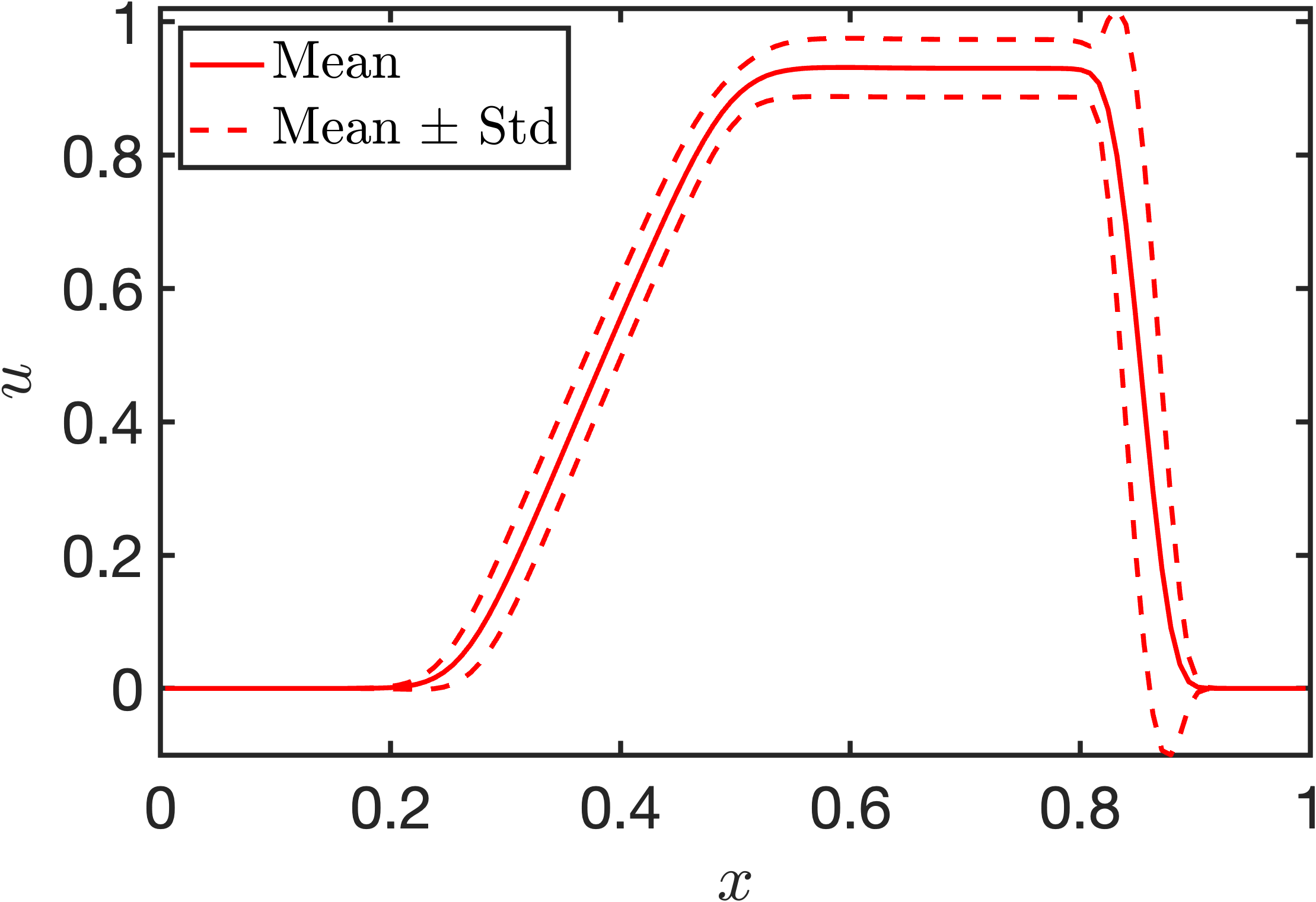}
        \caption{Mean$\pm$std of $u$}
    \end{subfigure}
    \caption{Stochastic Sod shock tube 2: solution plots and mean for full SFV method using WENO with reconstructed fluxes.}
    \label{fig:Sod2_flux}
\end{figure}

\begin{figure}
    \centering
    \begin{subfigure}{0.49\textwidth}
        \centering
        \includegraphics[width=\textwidth]{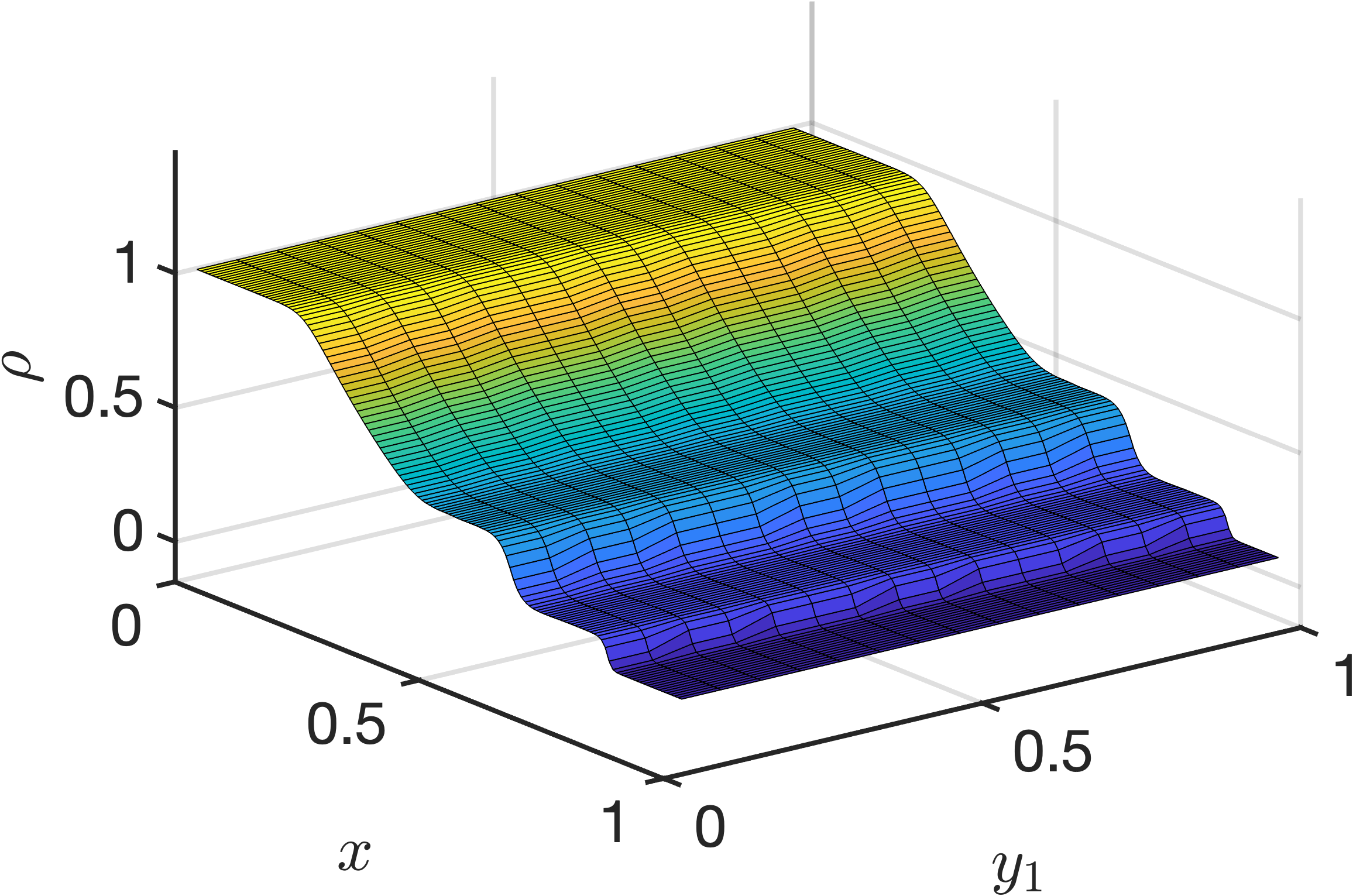}
        \caption{$\rho$ plot at $y_2 =1$, $N=30$ ROM}
    \end{subfigure}
    \begin{subfigure}{0.49\textwidth}
        \centering
        \includegraphics[width=\textwidth]{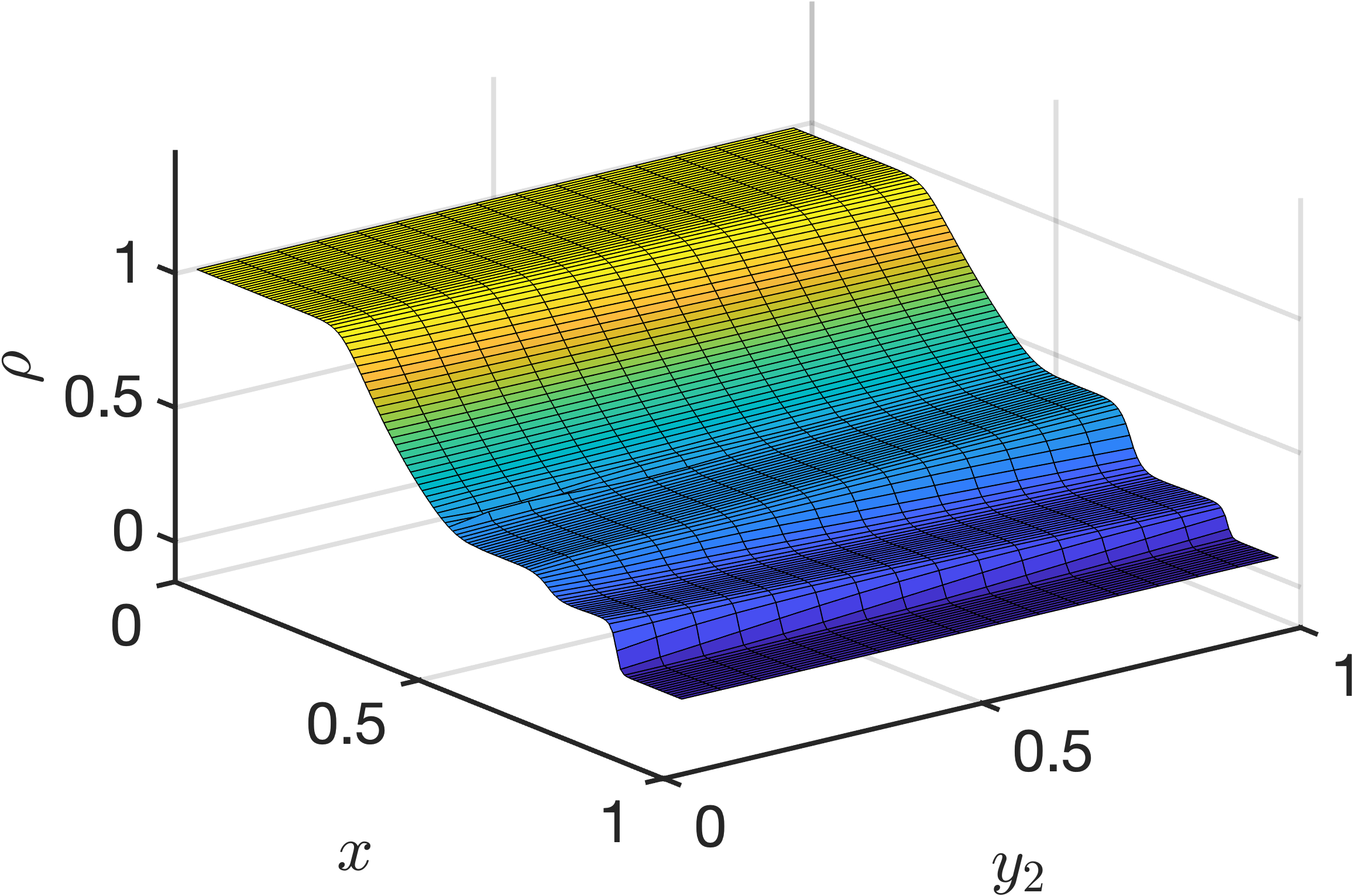}
        \caption{$\rho$ plot at $y_1 =1$, $N=30$ ROM}
    \end{subfigure}
     \begin{subfigure}{0.49\textwidth}
        \centering
        \includegraphics[width=\textwidth]{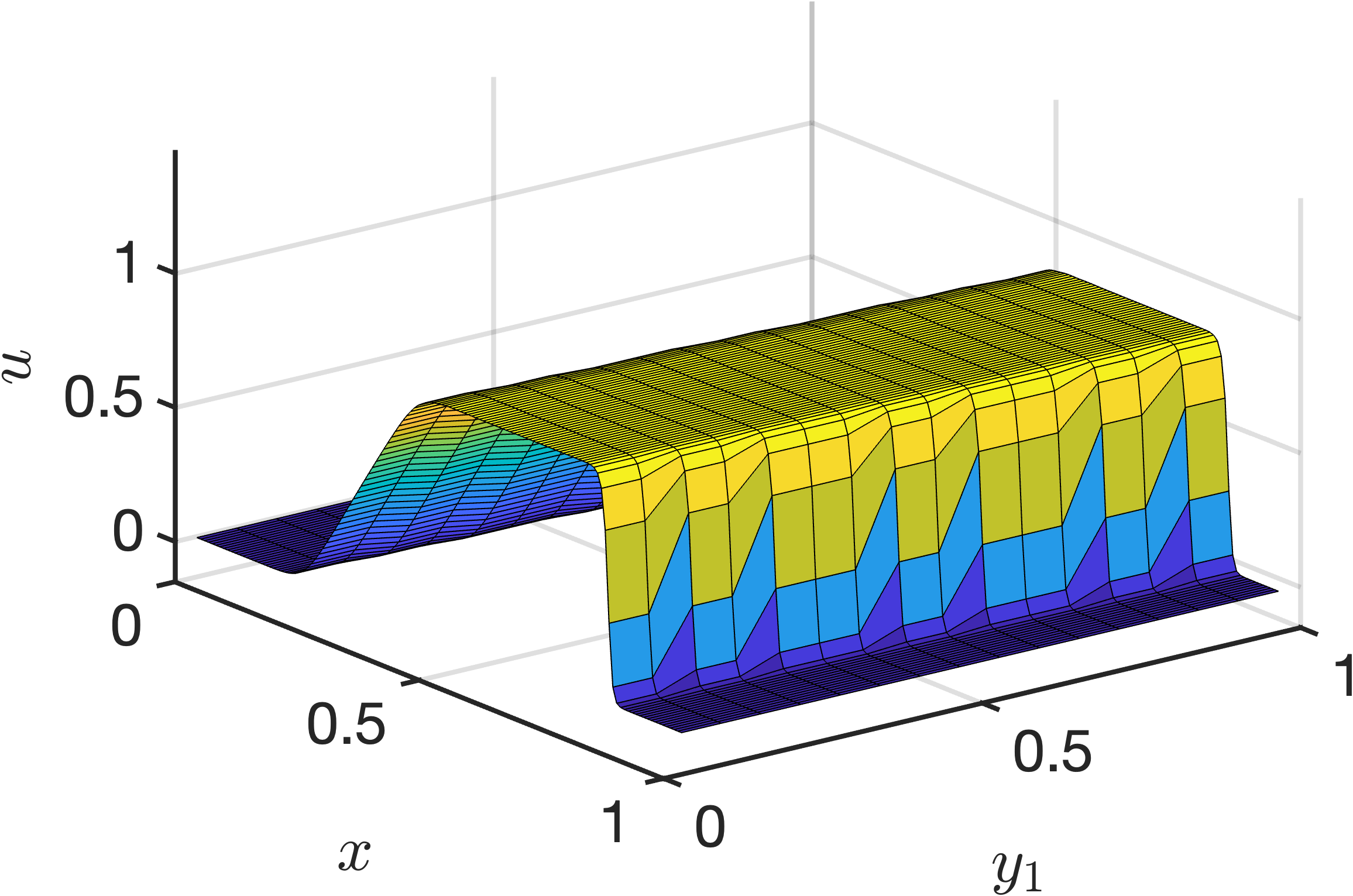}
        \caption{$u$ plot at $y_2 =1$, $N=30$ ROM}
    \end{subfigure}
    \begin{subfigure}{0.49\textwidth}
        \centering
        \includegraphics[width=\textwidth]{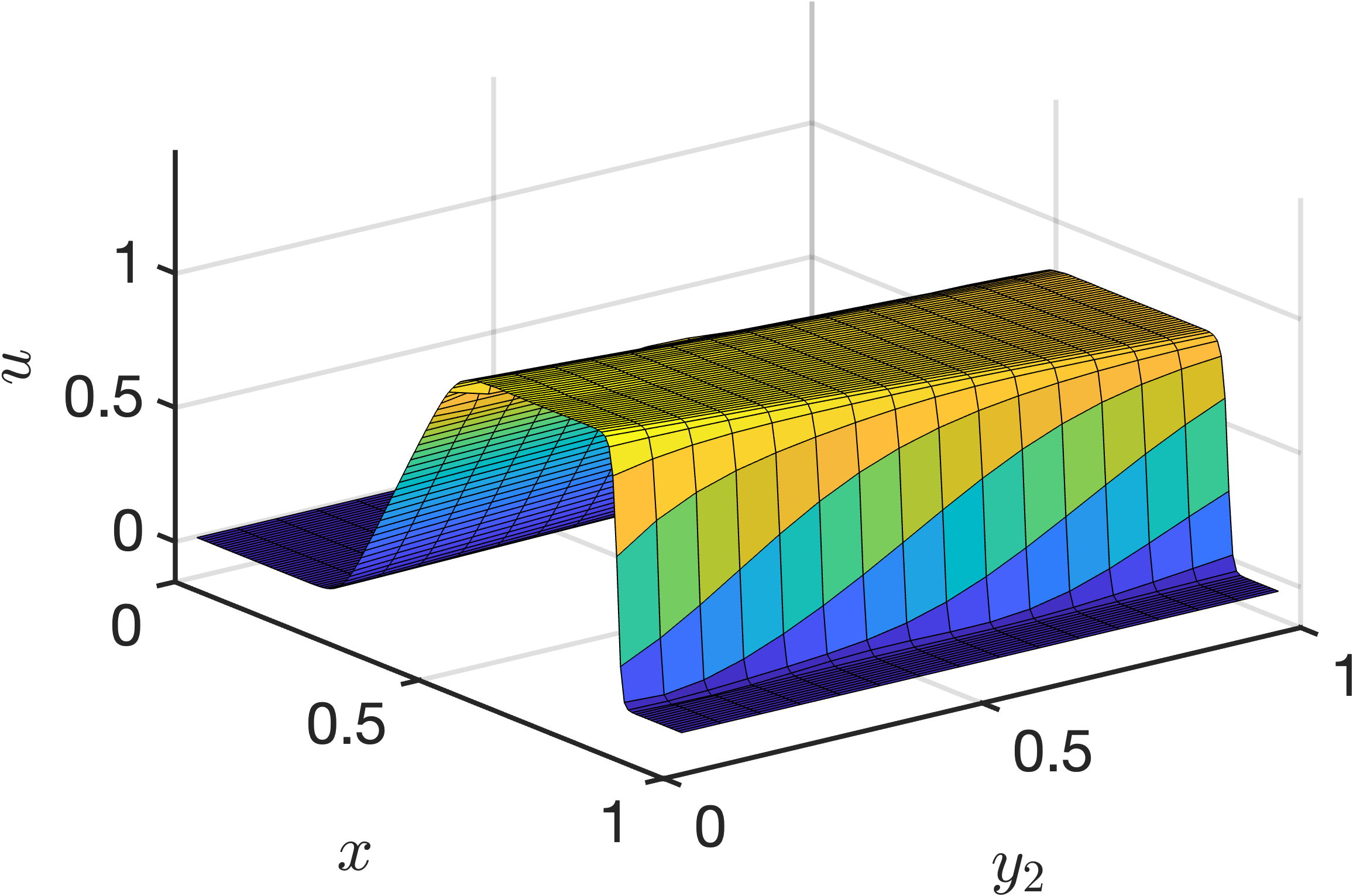}
        \caption{$u$ plot at $y_1 =1$, $N=30$ ROM}
    \end{subfigure}
     \begin{subfigure}{0.49\textwidth}
        \centering
        \includegraphics[width=\textwidth]{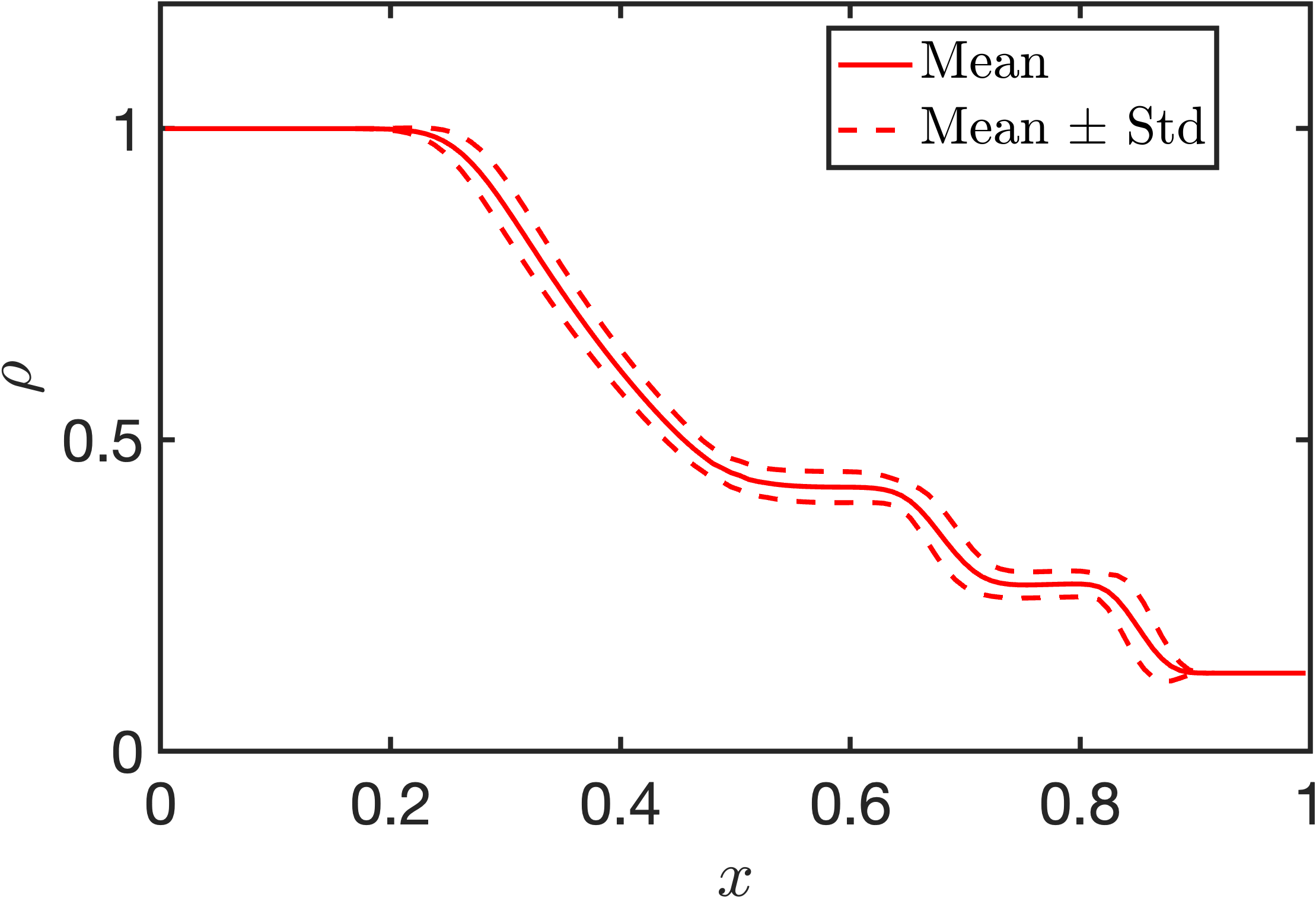}
        \caption{Mean$\pm$std of $\rho$, $N=30$ ROM}
    \end{subfigure}
    \begin{subfigure}{0.49\textwidth}
        \centering
        \includegraphics[width=\textwidth]{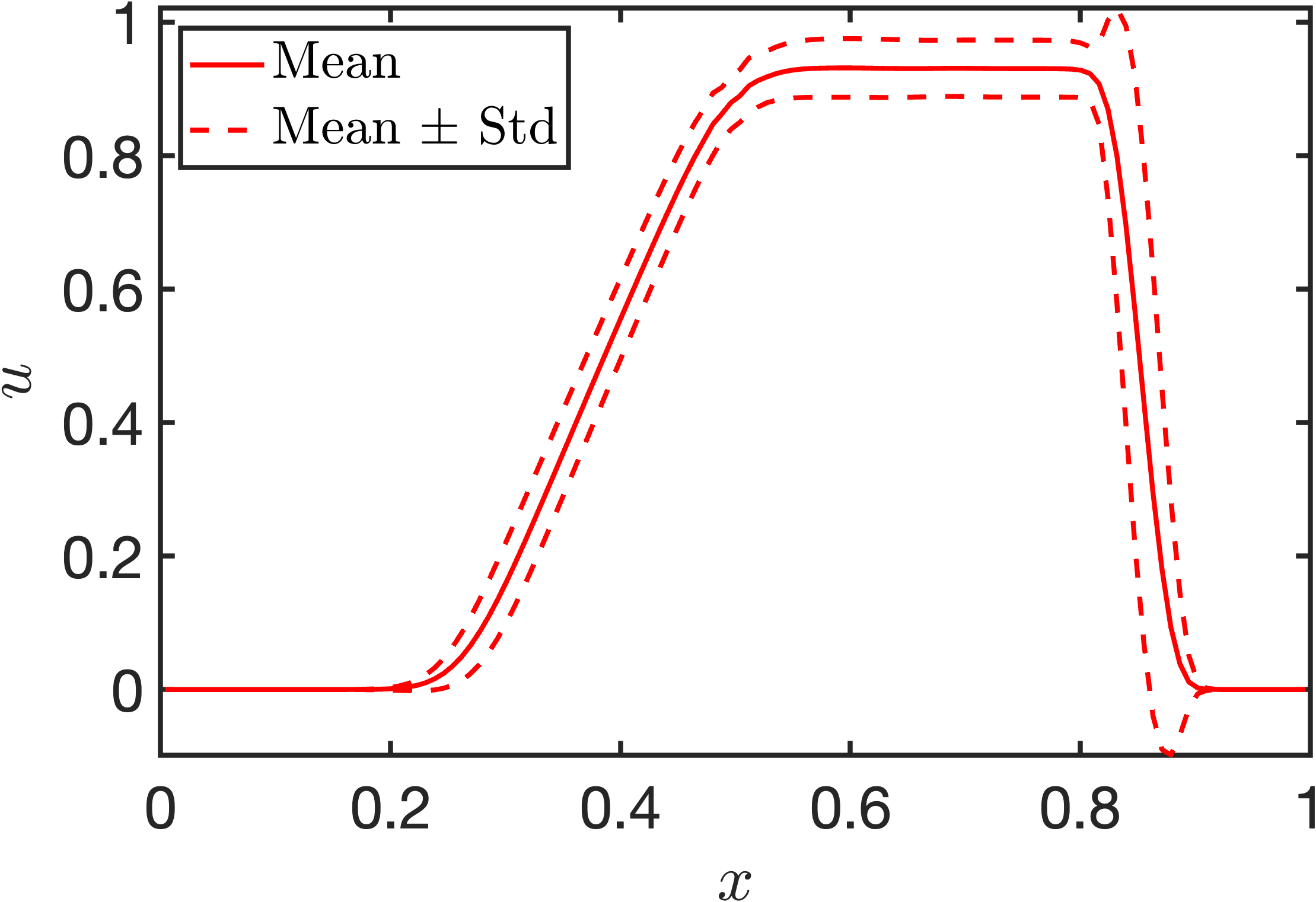}
        \caption{Mean$\pm$std of $u$, $N=30$ ROM}
    \end{subfigure}
    \caption{Stochastic Sod shock tube 2: solution plots and mean for $N_H = N = 30$ (resolved) hyper-reduced ROM.}
    \label{fig:Sod2_N30}
\end{figure}

\subsection{Sod shock tube with two stochastic variables}
We extend the stochastic Sod shock tube problem in Section~\ref{sec:Sod1} by introducing an additional stochastic variable to incorporate uncertainty in the flux functions. Specifically, in addition to the uncertainty in the initial shock position introduced through $y_1$ in \eqref{eq:Sod1_IC}, we model the adiabatic constant $\gamma$ as a random variable given by
\[
\gamma = y_2,
\]
where $y_2 \sim \mathcal{U}[1.2,1.6]$. Such randomness was considered for the SFV method in \cite{Tokareva14}.

Initially, we fix $N_x = 128$ and vary $N_y$ with the same discretization used for both stochastic variables $y_1$ and $y_2$ (each with size $\sqrt{N_y}$), and solve the problem using the full SFV method with both WENO reconstruction of states and WENO reconstruction of fluxes. We report their relative differences and the corresponding time-stepping runtimes in \autoref{tb:Sod2}. The results show that the two schemes produce close solutions with the relative difference decrease as $N_y$ increases, while WENO with reconstructed fluxes consistently achieves lower runtime.

Next, we fix $N_y = 16^2$ and examine the singular values of the flux snapshot matrix in \autoref{fig:Sod2_svd}. Compared to the previous single stochastic variable cases, the singular values exhibit a slower decay, with a noticeable drop occurring around $N = 200$. We then vary the number of modes $N$ from 15 to 40 in increments of 5, and report the hyper-reduced $N_H=N$ ROM error and time-stepping runtime in \autoref{fig:Sod2_er_time}. We observe a consistent decrease in error as $N$ increases. In contrast, the runtime first decreases and then increases, with ROMs in the range $20 \leq N \leq 30$ achieving runtime that is faster than that of the full SFV method.

Finally, we plot the solutions for $\rho$ and $u$ at $y_1 = 1$ and $y_2 = 1$ for the full SFV method using WENO with reconstructed fluxes, along with their corresponding mean and standard deviation, in \autoref{fig:Sod2_flux} as a reference for our ROM. We then generate the same set of plots for the ROM with $N=30$ in \autoref{fig:Sod2_N30}. At $N=30$, small oscillations are present in the solution plots, but the mean and standard deviation show no visible difference compared to the full SFV method.

\begin{remark}
   One might be concerned that the speedups achieved by the hyper-reduced ROM are not particularly substantial in our numerical experiments, and that for certain choices of modes $N$, the runtime per step can even be larger than that of the full SFV method due to insufficient modes to resolve shocks, causing the \texttt{ODE45} time-stepper to reject more steps. We emphasize, however, that the current implementation is not optimized. In particular, some routines used in our MATLAB implementation can be computationally inefficient, which further limits the observed performance gains.
\end{remark}

\section{Summary and future work}
In this work, we introduce an interpolation-based reduced order model for dimensionality reduction of the SFV method, designed to solve uncertainty quantification problems in systems of conservation laws. The model approximates flux integrals using a global stochastic POD basis, enabling efficient offline integration of basis functions while computing only the coefficients during the online stage. Additional computational savings are achieved through Q-DEIM hyper-reduction, which selects a subset of interpolation points from the full quadrature grid. We also propose a non-intrusive approach for constructing the snapshot matrix by solving multiple 1D problems using the standard FV method and applying WENO reconstructaion directly to the stochastic snapshots. Together, these techniques provide a scalable and accurate framework for uncertainty quantification.

In future work, we aim to enhance the accuracy of the ROM and extend its formulation to conservation laws in two- and three-dimensional physical domains. For example, while we use one set of reduced stochastic basis for all spatial points in this work, different sets of reduced stochastic bases could be computed independently at each spatial coordinate. This would adapt the dimension of the reduced stochastic basis to the stochastic solution at different spatial points, potentially improving accuracy. In addition, we plan to incorporate adaptive sampling strategies in the construction of stochastic reduced bases to improve the efficiency of the offline stage.

\section*{Declarations}

\subsection*{Funding}
Ray Qu and Jesse Chan gratefully acknowledge support from National Science Foundation under award DMS-1943186. Jesse Chan additionally acknowledges support from the National Science Foundation under award DMS-2231482.
Svetlana Tokareva gratefully acknowledges the support of the NNSA through the Laboratory Directed Research and Development (LDRD) program at Los Alamos National Laboratory under project number 20250032DR. Los Alamos National Laboratory is operated by Triad National Security, LLC for the U.S. Department of Energy’s NNSA. The Los Alamos unlimited release number is LA-UR-25-24954.

\subsection*{Availability of data and materials}
All simulation codes used in this work are publicly available at: \url{https://github.com/rayqu1126/SFV_ROM}. The code repository is also archived at: \url{https://doi.org/10.5281/zenodo.20078745}. Data supporting the results are available upon request.

\bibliographystyle{spmpsci}      % mathematics and physical sciences
\bibliography{sn-bibliography}   % name your BibTeX data base

\end{document}